\documentclass[11pt]{article}
\usepackage{amssymb,amsfonts,amsmath,amsthm,mathrsfs}
\usepackage{anysize}
\usepackage{hyperref}
\usepackage{color}
\usepackage{a4, fullpage}

\begin{document}
\baselineskip=14pt

\makeatletter\@addtoreset{equation}{section}
\makeatother\def\theequation{\thesection.\arabic{equation}}

\newcommand{\todo}[1]{{\tt [TODO: {#1}]}}

\newtheorem{defin}{Definition}[section]
\newtheorem{Prop}{Proposition}
\newtheorem{teo}{Theorem}[section]
\newtheorem{ml}{Main Lemma}
\newtheorem{con}{Conjecture}
\newtheorem{cond}{Condition}
\newtheorem{conj}{Conjecture}
\newtheorem{prop}[teo]{Proposition}
\newtheorem{lem}{Lemma}[section]
\newtheorem{rmk}[teo]{Remark}
\newtheorem{cor}{Corollary}[section]

\newcommand{\be}{\begin{equation}}
\newcommand{\ee}{\end{equation}}
\newcommand{\ben}{\begin{eqnarray}}
\newcommand{\benn}{\begin{eqnarray*}}
\newcommand{\een}{\end{eqnarray}}
\newcommand{\eenn}{\end{eqnarray*}}
\newcommand{\bp}{\begin{prop}}
\newcommand{\ep}{\end{prop}}
\newcommand{\bt}{\begin{teo}}
\newcommand{\et}{\end{teo}}
\newcommand{\bcor}{\begin{cor}}
\newcommand{\ecor}{\end{cor}}
\newcommand{\bcon}{\begin{con}}
\newcommand{\econ}{\end{con}}
\newcommand{\bcond}{\begin{cond}}
\newcommand{\econd}{\end{cond}}
\newcommand{\br}{\begin{rmk}}
\newcommand{\er}{\end{rmk}}
\newcommand{\bl}{\begin{lem}}
\newcommand{\el}{\end{lem}}
\newcommand{\bit}{\begin{itemize}}
\newcommand{\eit}{\end{itemize}}
\newcommand{\bd}{\begin{defin}}
\newcommand{\ed}{\end{defin}}
\newcommand{\bpr}{\begin{proof}}
\newcommand{\epr}{\end{proof}}

\newcommand{\Z}{{\mathbb Z}}
\newcommand{\R}{{\mathbb R}}
\newcommand{\E}{{\mathbb E}}
\newcommand{\C}{{\mathbb C}}
\renewcommand{\P}{{\mathbb P}}
\newcommand{\N}{{\mathbb N}}

\newcommand{\Bi}{{\cal B}}
\newcommand{\Si}{{\cal S}}
\newcommand{\Ti}{{\cal T}}
\newcommand{\Wi}{{\cal W}}
\newcommand{\Yi}{{\cal Y}}
\newcommand{\Hi}{{\cal H}}
\newcommand{\Fi}{{\cal F}}
\newcommand{\Zi}{{\cal Z}}

\newcommand{\eps}{\epsilon}

\newcommand{\nn}{\nonumber}

%Differentiation
\newcommand{\pa}{\partial}
\newcommand{\ffrac}[2]{{\textstyle\frac{{#1}}{{#2}}}}
\newcommand{\dif}[1]{\ffrac{\partial}{\partial{#1}}}
\newcommand{\diff}[1]{\ffrac{\partial^2}{{\partial{#1}}^2}}
\newcommand{\difif}[2]{\ffrac{\partial^2}{\partial{#1}\partial{#2}}}

%limits
\newcommand{\asto}[1]{\underset{{#1}\to\infty}{\longrightarrow}}
\newcommand{\Asto}[1]{\underset{{#1}\to\infty}{\Longrightarrow}}
\newcommand{\astoo}[1]{\underset{{#1}\to 0}{\longrightarrow}}
\newcommand{\Astoo}[1]{\underset{{#1}\to 0}{\Longrightarrow}}

\title{One-dimensional random walks 
with self-blocking immigration}
\author{Matthias Birkner$^{\,1}$, Rongfeng Sun$^{\,2}$}
\date{September 11, 2015}
\maketitle

\footnotetext[1]{Institut f\"ur Mathematik, Johannes-Gutenberg-Universit\"at Mainz, Staudingerweg 9, 55099 Mainz, Germany. Email:
birkner@mathematik.uni-mainz.de}

\footnotetext[2]{Department of Mathematics, National University of Singapore, 10 Lower Kent Ridge Road, 119076 Singapore. Email:
matsr@nus.edu.sg}

\begin{abstract}
We consider a system of independent one-dimensional random walkers 
where new particles are added at the origin at fixed rate whenever there 
is no older particle present at the origin. A Poisson ansatz leads to 
a semi-linear lattice heat equation and predicts that 
starting from the empty configuration the total number of 
particles grows as $c \sqrt{t} \log t$. We confirm this prediction 
and also describe the asymptotic macroscopic profile of the 
particle configuration. 
\end{abstract}

\section{Introduction: model and results}
%We are interested in 
Consider the following model of {\em random walks with self-blocking
immigration} (RWSBI) at the origin. Let $\eta_x(t)$ be the number of particles at position $x\in \Z$ at
time $t\geq 0$.  Particles perform independent continuous-time random walks on $\Z$ with jump rate $1$ and jump increments following a probability kernel
$(a_x)_{x\in\Z}$ with
\begin{equation}\label{kernela}
\sum_x x a_x = 0 \quad \mbox{and} \quad \sigma^2 := \sum_x x^2 a_x \in (0,\infty).
\end{equation}
In addition, at rate $\gamma > 0$ new particles attempt to ``immigrate'' at the origin $0$ but are
only successful if there is currently no other particle at $0$.  The
system starts with no particles at time $0$, i.e.\ $\eta_x(0) \equiv 0$. 
See Remark~\ref{rem:discuss} below for a discussion of the formal 
construction. 
\smallskip

This system shows interesting self-organized behavior: It possesses an
intrinsically defined ``correct'' growth rate and when particles are
added to the system at a lower (resp.\ higher) rate than this
correct rate, there is more (resp.\ less) vacant time at the
origin, which results in more (resp.\ less) particles added, and the
system is thus driven back toward the correct rate of addition
of particles. The task is thus to identify this correct asymptotic
rate at which particles are added to the system.
\smallskip

Obviously, more and more particles will be added to the system as time 
progresses and once created these perform independent random walks, which 
suggests hydrodynamic limit type arguments and results. While hydrodynamic 
limits for interacting particle systems is a vigorous area of 
current research, it seems that our system is somewhat special 
in this framework, and that there is presently no
readily applicable general theory to analyse its long-term behaviour: 
It combines a ``Kawasaki type'' dynamics, 
namely the motion of particles which preserves total mass, and
a very localised ``Glauber type'' dynamics, namely the immigration mechanism 
which creates new mass, in a non-trivial and non-perturbative way.
There is recent interest in extending hydrodynamic limits to models 
where non-trivial interactions among particles occur only in a very 
small part of the space, for example \cite{CF13} study systems of walks 
in bounded domains where pairwise annihilation only happens at the boundary. 
Insofar, our analysis of RWSBI fits to these efforts though our approach 
and the model details are quite different from \cite{CF13}. 
Arguably, RWSBI is of a very special form, yet 
we believe that at this stage, with no general approach available, 
a detailed analysis of special cases is warranted. 
\smallskip

Finally, we note that RWSBI first appeared in the literature 
as a caricature system for the analysis of a
system of critically branching random walks with a
density-dependent feedback, cf Remark~\ref{rmk:scbrw-rel} below.
\medskip

It is well known, see e.g.\ \cite[Ch.~1]{KL99}, that equilibrium states
for systems of independent random walks are products of Poisson distributions.
A Poisson ansatz leads to the heuristics that the
particle density $\E[\eta_x(t)] \approx \rho_x(t)$, where $\rho_x(t)$ is the
unique solution of the following ODE system, a semilinear discrete 
heat equation 
(the form of the non-linearity in the first line of \eqref{eqnrho}
arises by assuming $\eta_0(t)$ to be Poisson distributed with mean $\E[\eta_0(t)]$):
\begin{align}
\label{eqnrho}
\partial_t \rho_x(t) & = L_{\rm rw} \rho_x(t)
        +\gamma \delta_0(x) \exp(-\rho_0(t)),
        \quad t \geq 0, \: x \in \Z, \\ \notag
\rho_x(0) & \equiv  0, \quad x \in \Z,
\end{align}
where $L_{\rm rw}$ is the adjoint of the generator of the random walk given in 
(\ref{kernela}), with $(L_{\mathrm{rw}} f)_x := \sum_{y} a_{x-y} \big( f_y - f_x \big)$.

Denote the total mass of $\rho_\cdot(t)$ by
\begin{equation}\label{Rt}
R(t) := \sum_{x\in\Z} \rho_x(t) = \int_0^t \gamma \exp(-\rho_0(s)) \, ds.
\end{equation}
We have for $t\to\infty$
(see \cite[Lemma~17]{B03} and also Lemma~\ref{lemmarhoasympt} in
Appendix~\ref{Sect:rhoasymptotics})
\begin{eqnarray}
\label{rhod1localgrowth}
\rho_0(t) & = &
  \frac{1}{2}\log t - \log\log t
  + \log\big(\sqrt{2 \pi} \gamma /\sigma\big) + o(1), \\
\label{rhod1totalgrowth}
R(t) & = &  \Big[ \sigma\big(\tfrac{2}{\pi}\big)^{1/2} \sqrt{t} \log t
\Big] \big(1+ o(1) \big).
\end{eqnarray}
Furthermore (cf Lemma~\ref{lem:rhoscaling} below),
\begin{equation}
\label{eq:rhoscaling}
\frac1{\log t}\rho_{[\sigma\sqrt{t} y]}(t) \mathop{\longrightarrow}_{t\to\infty}
\tilde{\rho}(y) :=
\frac1{2\pi} \int_0^1 \frac1{\sqrt{s(1-s)}} e^{-y^2/(2s)}\, ds
= \frac1{\sqrt{2\pi}} \int_{|y|}^\infty e^{-z^2/2}\, dz,
\quad y \in \R.
\end{equation}
\smallskip

Our main result is that the Poisson ansatz is indeed valid. The asymptotic behavior of the total number of particles in the system, as well as the particle distribution in space, agree with the behavior of $\rho_\cdot(t)$ under the Poisson ansatz.

\begin{teo}
\label{thm:growth}
Let the model of random walks on $\Z$ with self-blocking immigration at the origin be defined as above, and recall $R(t)$ from \eqref{Rt} and \eqref{rhod1totalgrowth}. Almost surely, the total number of particles in the system satisfies
\begin{align}
\label{eq:conjgrowth}
\lim_{t\to\infty} \frac{1}{R(t)}\sum_x \eta_x(t) = 1.
\end{align}
\end{teo}

{Using Theorem~\ref{thm:growth}, we can further show that the ``shape of the particle cloud'', $(\eta_x(t))_{x\in \Z}$, follows the prediction from the Poisson ansatz.}
\begin{teo} \label{thm:shape} For any non-negative bounded continuous function $f \in C_{b,+}(\R)$, a.s.\ we have
\begin{align}
\label{eq:conjshape}
\lim_{t\to\infty} \frac{1}{\sigma\sqrt{t} \log t} \sum_x \eta_{x}(t) f\Big(\frac{x}{\sigma\sqrt t}\Big)  =
\int_\R f(y) \tilde \rho(y) \, dy,
\end{align}
where $\tilde \rho(y)= \frac1{\sqrt{2\pi}} \int_{|y|}^\infty e^{-z^2/2}\, dz$, as in \eqref{eq:rhoscaling}.
\end{teo}

\begin{rmk} \label{rem:discuss} \rm 
1.\ Starting from any finite initial condition, it is straightforward to construct 
the system $\eta$ explicitly by using suitable Poisson processes, for 
example as in Section~\ref{Sect:Poissyst} below; note that the total 
number of immigrated particles up to time $t$ is dominated by a rate $\gamma$ 
Poisson process, in particular the total number of particles is  a.s.\ finite
uniformly in any bounded time interval.

For a formal definition and suitable state space that allows infinite
configurations see \cite[Sect.~3.1]{B03}, compare also the arguments
in \cite[Sect.~2.2]{B03} for the construction of the transition
semigroup and a representation of $\eta$ as a Poisson process-driven
SDE system (a similar construction appears in~\cite{GK06}).
\smallskip

\noindent 
2.\ A much weaker version of \eqref{eq:conjgrowth} was previously shown in \cite[Prop.~8]{B03} via the relative entropy method \cite{Y91}, namely
that for any $\epsilon>0$,
\begin{equation}
\sum_{x\in\Z} \eta_x(t) = o(t^{1/2+\epsilon}) \quad
\text{in probability as $t\to\infty$}.
\end{equation}
%\smallskip

\noindent 
3.\ For the analogous system consisting of symmetric simple random
walks on $\Z^2$, a Poisson ansatz predicts $\rho_0(t) = \log\log t -
\log\log\log t - \log(2\pi) + o(1)$ and $R(t) \sim (2\pi t \log \log
t)/\log t$, cf \cite[Rem.~13]{B03}.  Using the techniques from
Section~\ref{Sect:ub}, it is fairly straightforward to establish
a corresponding upper bound for the total number of particles in the 
two-dimensional system in probability. 
It appears that in order to strengthen this bound to control the a.s.\ 
behavior and also to provide a matching lower bound using arguments parallel to 
those from Section~\ref{Sect:lb}, a very detailed study of the vacant time 
fluctuations of suitably tuned Poisson systems of two-dimensional random walks
with immigration will be required. We defer this question to future research. 
\end{rmk}

\begin{rmk}[Relation to self-catalytic branching random walks, 
\protect{\cite[Ch.~2]{B03}}] 
\label{rmk:scbrw-rel}
\rm 
Let SCBRW($b$) be a system of self-catalytic 
critical binary branching random walks on $\Z^d$ where each 
particle independently performs a random walk with kernel \eqref{kernela} 
and in addition while there are $k-1$ other particles at its site, 
it splits in two or disappers with rate $b(k)$, where $b : \N_0 \to 
[0,\infty)$ is a branching rate function (when $b$ is a linear function, 
this is a classical system of independent branching random walks).
Starting from a homogeneous initial condition, say a 
Poisson field on $\Z^d$ with constant intensity,
the long-term behaviour of such systems exhibits a dichotomy between 
persistence (i.e., convergence to a non-trivial shift-invariant equilibrium) 
and clustering (i.e., local extinction combined with increasingly rare regions 
of diverging particle density), depending on the branching rate function $b$ 
and the spatial dimension $d$. For general $b$ and $d\le 2$, it is 
believed (\cite[Conj.~1]{B03}) but not rigorously known that clustering occurs.
It is known, see e.g.\ \cite[Lemma~8]{B03}, that in this case clustering is 
equivalent to the local divergence as time $t\to\infty$ of the configuration 
under the so-called Palm distribution (which re-weights 
configurations at time $t$ proportional to the number of particles at 
the origin).

By a comparison result for the semigroups of SCBRW($b$) 
with respect to convex order for different $b$'s, cf 
\cite[Thm.~1 and Cor.~1]{B03}, it suffices to consider the special case 
$b=b_{\rm sing}$ with
$b_{\rm sing}(k)=\mathbf{1}_{\{k=1\}}$, i.e., particles branch
only if there is no other particle present at their site. 
The Palm distribution of SCBRW($b_{\rm sing}$) has a stochastic representation, 
\cite[Prop.~5]{B03}: It consists of the original SCBRW($b_{\rm sing}$) plus 
one special space-time path, which itself is drawn from the law of 
the time-reversed random walk, along which new particles immigrate at rate 
$1$ but only when there is no other, older particle already present 
at this site; the special path and the immigrating particles have an 
interpretation as the family decomposition for a focal particle picked 
at the origin at time $t$. 
While this is conceptually appealing, it appears currently still too complex 
to allow a rigorous analysis of its long-time behaviour. 

Thus, we consider the following simplification or caricature,
originally proposed by Anton Wakolbinger: Replace the random walk
special path by a constant path and disallow branching away from the
special path but keep the immigration mechanism along it 
unchanged. This yields RWSBI, our present object of study.  In this
sense, Theorems~\ref{thm:growth} and \ref{thm:shape} corroborate
Conjecture~1 from \cite{B03} in a quantitative way and in fact suggest
that the typical number of particles under the Palm distribution of
SCBRW($b_{\rm sing}$) should diverge like $\log t$ in $d=1$.  However,
undoing the caricature steps to convert our findings into an actual 
proof of this conjecture will require new arguments, which is currently
work in progress~\cite{BS15}.

\end{rmk}

\medskip

The rest of the paper is organized as follows. Section~\ref{Sect:Poissyst} introduces and analyzes
Poisson systems of random walks with immigration at the origin.
The upper (resp.\ lower) bounds in Theorems~\ref{thm:growth} and \ref{thm:shape} are proved in Section~\ref{Sect:ub} (resp.~Section~\ref{Sect:lb}) by suitable coupling
and comparison with the Poisson system of random walks. In Appendix~\ref{Sect:rhoasymptotics}, we derive the asymptotics\ \eqref{rhod1localgrowth}--\eqref{rhod1totalgrowth}, while in Appendix~\ref{Sect:Poiscor}, we derive an estimate for $k$-event
``correlation functions'' for Poisson processes.
\medskip

\section{Poisson systems of random walks} \label{Sect:Poissyst}
The key tool in our proof is an auxiliary Poisson system of random walks,
$\widetilde{\eta}=(\widetilde{\eta}_x(t))_{x\in\Z, t\ge 0}$, where particles immigrate at $x=0$ at time-dependent rate
$\beta(t)$, for some suitable $\beta : [0,\infty) \to (0,\infty)$. Once arrived, they follow independent continuous-time random walks
with jump rate $1$. By coupling such a Poisson system with random walks with self-blocking immigration (RWSBI), in particular, by coupling the times when the origin is vacant in each process, we can obtain bounds on the number of particles added to the RWSBI in terms of the Poisson system. We will choose $\beta(t)$ to be perturbations of the rate $\gamma e^{-\rho_0(t)}$ dictated by the Poisson ansatz in (\ref{eqnrho}).

We note that the system of random walks $\widetilde \eta$ can be characterized as a Poisson point process $\Pi$ on the set
$\mathcal{S}$ of all c\`adl\`ag paths $\cup_{t\geq 0} \{ \zeta : [t,\infty) \to \Z \}$ (denote the starting time of $\zeta$ by
$\tau_\zeta$),
with intensity measure
\begin{align}
\label{eq:Pois.intens}
\nu(d\zeta) = \beta(\tau_\zeta)d\tau_\zeta
\P\big(X \in d\zeta(\cdot-\tau_\zeta)\big)
\end{align}
where $X=(X_t)_{t\geq0}$ is the rate $1$ continuous time random walk as specified in (\ref{kernela}), starting
at $X_0=0$. Then
\[
\widetilde{\eta}_x(t) = \Pi\Big( \big\{ \zeta : \zeta(t)=x \big\}\Big),
\quad x \in \Z^d, \, t \ge 0,
\]
in particular, $\widetilde{\eta}_x(t)$ is Poisson distributed with
mean $\int_0^t \beta(u) p_x(t-u)\, du$, and
\[
\P(\widetilde{\eta}_x(t)=0) =
\exp\Big[ - \int_0^t \beta(u) p_x(t-u)\, du \Big],
\]
where
\begin{equation}
p_x(s): =\P(X_s=x).
\end{equation}
Apart from the number of particles added to the system by time $t$, we will also be interested in the amount of time at which the origin is vacant, i.e.,
\[
\widetilde{V}_{s,t} := \int_s^t \mathbf{1}_{\{\widetilde{\eta}_0(r)=0\}}\, dr, \quad 0 \le s \le t.
\]

We collect below results on the Poisson systems of random walks which we will need later. To prove the upper (resp.\ lower) bound in Theorems~\ref{thm:growth} and \ref{thm:shape}, it turns out that the appropriate choice of immigration rate $\beta(t)$ for the Poisson system $\widetilde \eta$ is
\begin{equation}\label{betavareps}
\beta^{(+\varepsilon)}(t) := (1+\varepsilon) \gamma e^{-\rho_0(t)}, \quad \mbox{resp.} \quad \beta^{(-\varepsilon)}(t) := (1-\varepsilon) \gamma e^{-\rho_0(t)},
\end{equation}
where $\varepsilon>0$ is chosen sufficiently close to $0$, and $\rho_0(t)$ is as in \eqref{eqnrho}. We will let $\widetilde \eta^{(\pm\varepsilon)}$ denote the respective Poisson system, and $\widetilde{V}_{s,t}^{(\pm\varepsilon)}$ its vacant time at the origin.

\begin{lem}\label{L:uppPoissonbd}
Let $\widetilde \eta^{(\pm\varepsilon)}$ be the Poisson system of random walks with immigration rate $\beta^{(\pm\varepsilon)}$ for some $\varepsilon\in (0,1)$. Then
\begin{align}
\label{eq:etaplus.totalno}
\frac{\sum_x \widetilde{\eta}^{(\pm\varepsilon)}_x(t)}{R(t)} \to 1\pm\varepsilon
\quad \text{a.s.\ as $t\to\infty$},
\end{align}
where $R_t=\sum_x \rho_x(t) =\big[ \sigma\big(\tfrac{2}{\pi}\big)^{1/2} \sqrt{t} \log t \big] \big(1+ o(1) \big)$ as defined in \eqref{Rt}, and
\begin{align}
\label{eq:etaplus.vacant}
\frac{\widetilde{V}^{(+\varepsilon)}_{0,t}}{R(t)} \to 0
\quad \text{a.s.\ as $t\to\infty$}.
\end{align}
\end{lem}

\begin{lem}\label{lem:lbPoiss.centr.k.mom}
Let $\widetilde \eta^{(-\varepsilon)}$ be the Poisson system of random walks with immigration rate $\beta^{(-\varepsilon)}$ for some $\varepsilon\in (0,1)$. Then there exists $t_0>0, c>0$ such that for all $t/2 \leq s < t$ with $t \geq t_0$, we have
\begin{align}
\label{eq:Evac.time.lb}
\E\big[\widetilde{V}^{(-\varepsilon)}_{s,t}] \geq c (t-s) t^{-\frac{1-\varepsilon}{2}}.
\end{align}
If $\xi \in (\tfrac{1-\varepsilon}{2},1)$, then there exists $b \in (0,\infty)$ such that for every
$k \in \N$,
there exist $t_0, C \in (0,\infty)$ so that for all $s,t$ with
$t_0 \leq t/2 \leq s \leq t-t^\xi$,
\begin{align} \label{centeredmombd}
  \E\Big[ \big( \widetilde{V}^{(-\varepsilon)}_{s,t}
    - \E\big[\widetilde{V}^{(-\varepsilon)}_{s,t}] \big)^k \Big]
  \leq C t^{-b k} \E\big[\widetilde{V}^{(-\varepsilon)}_{s,t}]^k
\end{align}
$($we can choose $b=\left(\xi-\tfrac{1-\varepsilon}{2} \right)/48\,)$.
\end{lem}
\noindent
This shows that the vacant time $\widetilde{V}^{(-\varepsilon)}_{s,t}$
is concentrated around its mean with high probability.
\medskip

We now give the proofs of Lemmas~\ref{L:uppPoissonbd} and \ref{lem:lbPoiss.centr.k.mom}. The proof of Lemma~\ref{lem:lbPoiss.centr.k.mom}
is of independent interest, but is quite involved, and therefore can be read after the proof of Theorems~\ref{thm:growth} and \ref{thm:shape} in Sections~\ref{Sect:ub} and \ref{Sect:lb}.

\bigskip
\noindent
{\bf Proof of Lemma~\ref{L:uppPoissonbd}.}
Recalling \eqref{betavareps}, we have
$$
\widetilde{\rho}^{(\pm\varepsilon)}_x(t) :=
\E[\widetilde{\eta}^{(\pm\varepsilon)}_x(t)]
= \int_0^t \beta^{(\pm\varepsilon)}(s) p_x(t-s) \, ds
= (1\pm\varepsilon) \int_0^t \gamma e^{-\rho_0(t)} p_x(t-s) \, ds
= (1\pm\varepsilon) \rho_x(t).
$$
In particular $\E\big[ \sum_x \widetilde{\eta}^{(\pm\varepsilon)}_x(t) \big]
= (1\pm\varepsilon)\sum_x \rho_x(t) = (1\pm\varepsilon) R(t)$. Since $\sum_x \widetilde{\eta}^{(\pm\varepsilon)}_x(t)$
is nothing but a time-changed Poisson process with mean $(1\pm\varepsilon)R(t)$, (\ref{eq:etaplus.totalno}) follows
immediately.
\smallskip

To prove (\ref{eq:etaplus.vacant}), note that
$\beta^{(\pm\varepsilon)}(t) \sim \frac{1\pm\varepsilon}{\sqrt{2\pi}} \sigma t^{-1/2} \log t$ by \eqref{rhod1localgrowth}, and hence
$$
\E\big[ \widetilde{V}^{(+\varepsilon)}_{0,t} \big] =
\int_0^t e^{-\widetilde{\rho}^{(+\varepsilon)}_0(u)} \, du
= \int_0^t e^{-(1+\varepsilon) \rho_0(u)} \, du
\leq 1 + C \int_1^t \frac{(\log u)^{1+\varepsilon}}{u^{(1+\varepsilon)/2}} \, du \leq
2C t^{(1-\varepsilon)/2} (\log t)^{1+\varepsilon}.
$$
For any $\delta>0$, by Markov inequality and the asymptotics of $R(t)$ in (\ref{rhod1totalgrowth}),
$$
\P(\widetilde{V}^{(+\varepsilon)}_{0,t} > \delta R(t)) \leq \frac{\E\big[ \widetilde{V}^{(+\varepsilon)}_{0,t} \big]}{\delta R(t)}
\leq \frac{C'}{t^{\varepsilon/4}}
$$
for all $t$ sufficiently large. By Borel-Cantelli, along the sequence of times $t_n=c^n$ for any $c>1$, we then have $\limsup_{n\to\infty} \widetilde{V}^{(+\varepsilon)}_{0,c^n}/R(c^n) \leq \delta$ almost surely. Since $\widetilde{V}^{(+\varepsilon)}_{0,s} \le \widetilde{V}^{(+\varepsilon)}_{0,t}$ for $s \le t$, together with the asymptotics of $R(t)$ given in (\ref{rhod1totalgrowth}), we obtain
\[
\limsup_{t\to\infty} \frac{\widetilde{V}^{(+\varepsilon)}_{0,t}}{R(t)}
\le \limsup_{n\to\infty} \frac{\widetilde{V}^{(+\varepsilon)}_{0,c^n}}{R(c^{n-1})}
\leq \delta \sqrt{c}.
\]
Since $\delta>0$ can be chosen arbitrarily, (\ref{eq:etaplus.vacant}) then follows.
\qed
\bigskip

\noindent
{\bf Proof of Lemma~\ref{lem:lbPoiss.centr.k.mom}.} Using the asymptotics of $\rho_0(\cdot)$ given in (\ref{rhod1localgrowth}), (\ref{eq:Evac.time.lb}) holds because
$$
\begin{aligned}
\E[\widetilde{V}^{(-\varepsilon)}_{s,t}] =
\int_s^t e^{-\E[\widetilde \eta^{(-\varepsilon)}_0(u)]} \, du & = \int_s^t e^{-(1-\varepsilon) \rho_0(u)} \, du \\
& \geq \int_s^t u^{-\frac{1-\varepsilon}{2}} \, du =
\frac{2}{1+\varepsilon}
\big( t^{\frac{1+\varepsilon}{2}} - s^{\frac{1+\varepsilon}{2}} \big)
\geq c (t-s) t^{-\frac{1-\varepsilon}{2}}.
\end{aligned}
$$

Next we prove the centered moment bound (\ref{centeredmombd}). To lighten notation, we will drop the dependence on $\varepsilon$
in the remainder of the proof and write $\widetilde{V}_{s,t}=\widetilde{V}_{s,t}^{(-\varepsilon)}$, $\widetilde{\eta}_0=\widetilde{\eta}_0^{(-\varepsilon)}$, etc. Note
\begin{align}
\label{eq:lbPoiss.centr.k.mom1}
  \E\Big[ \big( \widetilde{V}_{s,t} - \E\big[\widetilde{V}_{s,t}] \big)^k \Big]
 = k! \hspace{-1.1em} \mathop{\int\cdots\int}_{s \leq u_1 < \cdots < u_k \leq t}
 \E\Big[ \prod_{i=1}^k \big( \mathbf{1}(\widetilde{\eta}_0(u_i)=0)
- \P(\widetilde{\eta}_0(u_i)=0) \big) \Big] \, du_k \dots du_1.
\end{align}
The idea to estimate \eqref{eq:lbPoiss.centr.k.mom1} is the following.
When the $u_i$'s are close, the contribution to the integral is small due to the restricted
range of integration; when the $u_i$'s are far apart, we can use the decorrelation
of the Poisson system as quantified by Lemma~\ref{lem:Poisvac}.
We thus group $u_i$'s into blocks as follows, where each block contains consecutive $u_i$'s that are close to each other,
and different groups are far apart.
\smallskip

We group the time points $u_1,\ldots, u_k$ into blocks that are separated from each other by at least $t^\delta$, with
$\delta = \frac23 \left( \xi-\tfrac{1-\varepsilon}{2}\right) \, (>0)$. A block structure is determined by $\ell \in \{1,2,\dots,\lfloor k/2\rfloor\}$, and $\ell$ pairs of indices $g_i, h_i$ with
$$
1 \leq g_1 < h_1 < g_2 < h_2 < \cdots < g_\ell < h_\ell \leq k.
$$
Let $B(g,h)$ denote the set of all $\vec u:=(u_1,\dots,u_k)$ with $s \leq u_1 < \cdots < u_k \leq t$, such that for each $1\leq i\leq \ell$,
$J_i:=[g_i, h_i]\cap \N$ is a block of indices with
$$
u_{r+1}-u_r \leq t^\delta \mbox{ for all } r\in [g_i, h_i-1]\cap \N \quad \mbox{and} \quad \min\{u_{g_i}-u_{g_i-1}, u_{h_i+1}-u_{h_i}\}> t^\delta,
$$
where $u_0:=-\infty$, $u_{k+1}:=+\infty$. Indices in the set $J_0 := \{ 1,\dots,k\} \setminus \big(J_1 \cup \cdots \cup J_\ell\big)$ are the blocks of singletons, i.e., for each $i\in J_0$, $u_i$ is separated from all the other $u_j$'s by at least $t^\delta$.

Now consider a fixed block structure as determined by $\ell$ and $g_1, h_1,\ldots, g_\ell, h_\ell$, and let $(u_1, \ldots, u_k)\in B(g, h)$.
Write
\begin{align}
\label{eq:prod.i.to.k}
\prod_{i=1}^k \big( \mathbf{1}(\widetilde{\eta}_0(u_i)=0)
- \P(\widetilde{\eta}_0(u_i)=0) \big)
= \prod_{m=0}^\ell
\prod_{i \in J_m} \big( \mathbf{1}(\widetilde{\eta}_0(u_i)=0)
- \P(\widetilde{\eta}_0(u_i)=0) \big).
\end{align}
To apply Lemma~\ref{lem:Poisvac}, for each block $J_m$ with $1\leq m\leq \ell$, we need to rewrite the product of the centered indicators as linear combinations of centered indicators. More precisely, for each $1\leq m\leq \ell$, we write
\begin{align}
& \prod_{i \in J_m} \big( \mathbf{1}(\widetilde{\eta}_0(u_i)=0)
- \P(\widetilde{\eta}_0(u_i)=0) \big) \notag \\
=\ & \sum_{J_m' \subset J_m} (-1)^{|J_m \setminus J_m'|}
\mathbf{1}(\widetilde{\eta}_0(u_i)=0, i \in J_m') \hspace{-0.6em}
\prod _{i \in J_m \setminus J_m'} \hspace{-0.8em} \P(\widetilde{\eta}_0(u_i)=0)
\notag \\
=\ & \sum_{J_m' \subset J_m} (-1)^{|J_m \setminus J_m'|}
\Big( \mathbf{1}(\widetilde{\eta}_0(u_i)=0, i \in J_m')
- \P\big( \widetilde{\eta}_0(u_i)=0, i \in J_m' \big) \Big)
\hspace{-0.6em}
\prod _{i \in J_m \setminus J_m'} \hspace{-0.8em} \P(\widetilde{\eta}_0(u_i)=0)
\notag \\
\label{eq:prod.in.Jm}
& \hspace{1.5em} {} +
\sum_{J_m'' \subset J_m} (-1)^{|J_m \setminus J_m''|}
\P\big( \widetilde{\eta}_0(u_i)=0, i \in J_m'' \big)
\hspace{-0.6em}
\prod _{i \in J_m \setminus J_m''} \hspace{-0.8em} \P(\widetilde{\eta}_0(u_i)=0),
\end{align}
where we centered the indicator function $\mathbf{1}(\widetilde{\eta}_0(u_i)=0, i \in J_m')$, and $\P( \widetilde{\eta}_0(u_i)=0, i \in J_m')$ is interpreted to be $1$ if $J'_m=\emptyset$. Note that for blocks of singletons, i.e., $i\in J_0$, the indicator function is already centered and there is no constant term as in (\ref{eq:prod.in.Jm}), which is why the singleton blocks are separated from the other blocks $J_1, \ldots, J_\ell$.

Applying \eqref{eq:prod.in.Jm} for indices in blocks $J_1, \ldots, J_\ell$ in \eqref{eq:prod.i.to.k}, and expanding and grouping terms,
we can then rewrite \eqref{eq:prod.i.to.k} as a sum of
\[
\pm A(\widehat{J},\vec J', \vec J'', \vec u),
\]
where $\widehat{J} \subset \{1,\dots,\ell\}$ determine the blocks for which we choose a centered indicator function (instead of a constant) from the expansion in (\ref{eq:prod.in.Jm}), and $\vec J'=(J'_1,\ldots, J'_\ell)$, $\vec J''=(J''_1,\ldots, J''_\ell)$, with $J_m', J_m'' \subset J_m$ as in (\ref{eq:prod.in.Jm}) for each block $J_m$. More precisely,
%\begin{equation}

%\begin{split}
\begin{align}
 & A(\widehat{J},\vec J', \vec J'', \vec u) \notag\\
=\ & \prod_{i \in J_0} \big( \mathbf{1}(\widetilde{\eta}_0(u_i)=0)
- \P(\widetilde{\eta}_0(u_i)=0) \big) \notag\\
& \times \prod_{m \in \widehat{J}} \bigg(
\Big( \mathbf{1}(\widetilde{\eta}_0(u_i)=0, i \in J_m')
- \P\big( \widetilde{\eta}_0(u_i)=0, i \in J_m' \big) \Big)
\hspace{-0.6em}
\prod_{i \in J_m \setminus J_m'} \hspace{-0.8em} \P(\widetilde{\eta}_0(u_i)=0)
\bigg) \notag\\
& \times \hspace{-1.2em}
\prod_{m \in \{1,\dots,\ell\} \setminus \widehat{J}} \bigg(
\P\big( \widetilde{\eta}_0(u_i)=0, i \in J_m'' \big)
\hspace{-0.6em}
\prod _{i \in J_m \setminus J_m''} \hspace{-0.8em} \P(\widetilde{\eta}_0(u_i)=0)
\bigg) \notag \\
=\ &
\prod_{i \in J_0} \big( \mathbf{1}(\widetilde{\eta}_0(u_i)=0)
- \P(\widetilde{\eta}_0(u_i)=0) \big)
\times \prod_{m \in \widehat{J}} \bigg(
\Big( \mathbf{1}(\widetilde{\eta}_0(u_i)=0, i \in J_m')
- \P\big( \widetilde{\eta}_0(u_i)=0, i \in J_m' \big) \Big) \notag \\
& \hspace{2em} \times \hspace{-1.2em}
\prod_{m \in \{1,\dots,\ell\} \setminus \widehat{J}}
\P\big( \widetilde{\eta}_0(u_i)=0, i \in J_m'' \big)
\: \times \:\prod_{i \in \breve{J}} \P(\widetilde{\eta}_0(u_i)=0), \label{eq:A.manyJ}
\end{align}
%\end{split}
%\end{equation}
where in the last line, $\breve{J} := \cup_{m \in \widehat{J}} \, (J_m \setminus J_m')
\, \bigcup \, \cup_{m \in \{1,\dots,\ell\} \setminus \widehat{J}} \,
(J_m \setminus J_m'')$.

The sign corresponding to a given choice of $\widehat{J} , \vec J', \vec J''$ is
\[
(-1)^{\sum_{m \in \widehat{J}} |J_m \setminus J_m'| +
\sum_{m \in \{1,\dots,\ell\} \setminus\widehat{J}} |J_m \setminus J_m''|} .
\]

Using Lemma~\ref{lem:Poisvac}, we can bound the expectation of the product of centered indicator functions in \eqref{eq:A.manyJ} by
\begin{align}
& \E\bigg[ \prod_{i\in J_0}
\Big( \mathbf{1}(\xi(E_i)=0) - \P(\xi(E_i)=0) \Big) \prod_{m\in \widehat J} \Big( \mathbf{1}(\xi(F_m)=0) - \P(\xi(F_m)=0) \Big) \bigg] \notag \\
\leq\ &
C t^{-\delta (|J_0|+|\widehat{J}|)/8}
\times \prod_{i \in J_0} \P(\widetilde{\eta}_0(u_i)=0)
\: \times \:
\prod_{m \in \widehat{J}} \P\big( \widetilde{\eta}_0(u_i)=0, i \in J_m' \big), \label{B1app1}
\end{align}
where $\xi$ is a Poisson point process on the random walk paths space with intensity measure $\nu$ given by
\eqref{eq:Pois.intens}, with $\beta(t)=\beta^{(-\varepsilon)}(t)$; and
\begin{align*}
E_i & := \{\text{random walk paths $\zeta$ : with}\; \zeta(u_r) = 0 \}, \quad r\in J_0,  \\
F_m & := \{\text{random walk paths $\zeta$ : with}\;
\zeta(u_r) = 0 \; \text{for some} \; r \in J_{m}' \}, \quad
m\in \widehat J.
\end{align*}
Let us reorder and relabel the sets $(E_i)_{i\in J_0}$ and $(F_m)_{m\in \widehat J}$ by $(\widetilde E_i)_{1\leq i\leq |J_0|+|\widehat J|}$, where each $\widetilde E_i$ is of the form $\{\zeta: \zeta(u_r)=0 \mbox{ for some } r\in \widetilde J_i\}$ for some distinct index set $\widetilde J_i\subset \{1, \ldots, k\}$, and elements of $\widetilde J_1$ being smaller than those of $\widetilde J_2$, etc.

To see how does (\ref{B1app1}) follow from Lemma~\ref{lem:Poisvac}, note that for any $I \subset \{1,2,\dots,|J_0|+|\widehat{J}|\}$,
$$
\begin{aligned}
\nu\big( \bigcap_{i \in I} \widetilde E_i\big) =\ & \nu(\zeta: \mbox{ for each } i\in I, \ \zeta(u_r)=0 \ \mbox{for some } r\in \widetilde J_i) \\
=\ & \int_0^t \beta^{(-\varepsilon)}(v) \P(X_{u_r}=0 \mbox{ for some } r\in \widetilde J_i \mbox{ for each } i\in I | X_v=0) dv \\
\leq \ & \Big(\frac{C'}{t^{\delta/2}}\Big)^{|I|-1}\int_0^t \beta^{(-\varepsilon)}(v) \P(X_{u_r}=0 \mbox{ for some } r\in \widetilde J_1 | X_v=0) dv \\
\leq \ &  \Big(\frac{C'}{t^{\delta/2}}\Big)^{|I|-1} \int_0^t \beta^{(-\varepsilon)}(v) \sum_{r\in \widetilde J_1}\P(X_{u_r}=0| X_v=0) dv \\
=\ & \Big(\frac{C'}{t^{\delta/2}}\Big)^{|I|-1} \sum_{r\in \widetilde J_1} \E[\widetilde \eta_0(u_r)] = (1-\varepsilon) \Big(\frac{C'}{t^{\delta/2}}\Big)^{|I|-1} \sum_{r\in \widetilde J_1} \rho_0(u_r) \leq \frac{C}{(t^{\delta/4})^{|I|-1}},
\end{aligned}
$$
where we applied the local central limit theorem in the first inequality, noting that the random walk returns to the origin at least $|I|-1$ times over intervals of length at least $t^\delta$, and we applied (\ref{rhod1localgrowth}) to bound $\rho_0(u_r)$ in the last inequality. For any $I_1, \dots, I_n \subset\{1,\dots,|J_0|+|\widehat{J}|\}$ with $|I_1|,\dots,|I_n| \geq 2$ and $I_1 \cup \cdots \cup I_n = \{1,2,\dots,|J_0|+|\widehat{J}|\}$, we then have
$$
\prod_{j=1}^n \nu\big( \bigcap_{i \in I_j} \widetilde E_i\big) \leq C^n \big(t^{-\delta/4}\big)^{\sum_{i=1}^n |I_i|-n} \leq \min\Big\{C^n (t^{-\delta/4})^{|J_0|+|\widehat J|-n}, (Ct^{-\delta/4})^n\Big\}.
$$
Substituting these bounds into (\ref{eq:Poisvac.full}), where the first bound is used for $1\leq n\leq \frac{|J_0|+|\widehat J|}{2}$, and the second bound is used for $n> \frac{|J_0|+|\widehat J|}{2}$, it is then easily seen that (\ref{B1app1}) follows (note that we only need to consider
(\ref{B1app1}) for the case $|J_0|+|\widehat J|\geq 2$, since otherwise the inequality is trivial).

Having verified (\ref{B1app1}), we can then apply (\ref{eq:A.manyJ}) to  bound
\begin{align}
& \E\Big[
A(\widehat{J}, \vec J', \vec J'', \vec u) \Big] \notag \notag \\
\leq\ & C t^{-\delta (|J_0|+|\widehat{J}|)/8}  \prod_{i \in J_0\cup \breve J} \P(\widetilde{\eta}_0(u_i)=0)
\prod_{m \in \widehat{J}} \P\big( \widetilde{\eta}_0(u_i)=0, i \in J_m' \big) \!\!\!\! \prod_{m \in \{1,\dots,\ell\} \setminus \widehat{J}}\!\!\!\!
\P\big( \widetilde{\eta}_0(u_i)=0, i \in J_m'' \big) \notag\\
\leq\ & C t^{-\delta |J_0|/8}  \prod_{i \in J_0\cup J_0'} \P(\widetilde{\eta}_0(u_i)=0)
\end{align}
where $J_0'$ contains the smallest index from each block $J_m$, $1\leq m\leq \ell$.

Therefore following the discussion after (\ref{eq:prod.in.Jm}), we have
\begin{align}
& \mathop{\int\cdots\int}_{B(g, h)}
 \E\Big[ \prod_{i=1}^k \big( \mathbf{1}(\widetilde{\eta}_0(u_i)=0)
- \P(\widetilde{\eta}_0(u_i)=0) \big) \Big] \, du_k \dots du_1
\notag \\
\leq\ &  C t^{-\delta |J_0|/8} \sum_{\widehat J, \vec J', \vec J''}
 \mathop{\int\cdots\int}_{B(g,h)}
\prod_{i \in J_0 \cup J_0'} \P(\widetilde{\eta}_0(u_i)=0) \, du_k \dots du_1
\notag \\
\leq\ &  C' t^{-\delta |J_0|/8} t^{\delta(|J_1|+\cdots +|J_\ell|-\ell)}
\E\big[ \widetilde{V}_{s,t} \big]^{|J_0|+\ell} \notag \\
=\ & C' t^{-\delta |J_0|/8} \left( \frac{t^\delta}{\E\big[ \widetilde{V}_{s,t} \big]}\right)^{k-|J_0|-\ell}
\E\big[ \widetilde{V}_{s,t} \big]^k, \label{Bghbd1}
\end{align}
where $C'$ contains combinatorial factors that depend only on $k$, but not on $s$ and $t$
(and we used $|J_1|+\cdots+|J_\ell|=k-|J_0|$ in the last line).

Since $\E\big[ \widetilde{V}_{s,t} \big] \geq
c (t-s) t^{-\frac{1-\varepsilon}{2}} \geq
c t^{\xi-\frac{1-\varepsilon}{2}}$ by \eqref{eq:Evac.time.lb}
and the assumption on $s$ and $t$,
the term in (\ref{Bghbd1}) is bounded by
\begin{align}
\label{eq:boundbyEtok.lb}
C'' t^{- \delta|J_0|/8} \ t^{(\delta-\xi+\frac{1-\varepsilon}{2}) (k-|J_0|-\ell)}
\ \E\big[ \widetilde{V}_{s,t} \big]^k .
\end{align}
Note that $\delta-\xi+\frac{1-\varepsilon}{2} < 0$, $k-|J_0|-\ell \geq 0$, and $\ell \leq (k-|J_0|)/2$ (since each block $J_m$, $1\leq m\leq \ell$, contains at least two indices). Thus when $|J_0| \geq k/4$, we can bound \eqref{eq:boundbyEtok.lb} by
\begin{align}
C'' t^{- k \delta/32} \E\big[ \widetilde{V}_{s,t} \big]^k ,
\end{align}
whereas when $|J_0| < k/4$ (and hence $k-|J_0|-\ell \geq k/4$), we can bound \eqref{eq:boundbyEtok.lb} by
\begin{align}
C'' t^{- k(\xi-\frac{1-\varepsilon}{2}-\delta)/4} \E\big[ \widetilde{V}_{s,t} \big]^k .
\end{align}
Either way, we find that the bound in (\ref{Bghbd1}) can be bounded by $C''t^{- bk} \E[ \widetilde{V}_{s,t}]^k$ for some $C''$ depending only $k$, and $b>0$ depending only on $\xi$ and $\varepsilon$.

Since for given $k$ there are only finitely many choices for $\ell$ and $g_1,h_1;\dots;g_\ell,h_\ell$,
summing over all possible $B(g,h)$ then yields the claimed bound (\ref{centeredmombd}) with $b=\delta/32 =\left(\xi-\tfrac{1-\varepsilon}{2} \right)/48$.
\qed

% %%%%%%%%%%%%%%%%%%%%%%%%%%% Upper bd %%%%%%%%%%%%%%%%%%%%%%%%%%

\section{Upper bounds in Theorems~\ref{thm:growth} and \ref{thm:shape}}
\label{Sect:ub}

Here is the basic idea for the upper bound on the system of random walks with self-blocking immigration (RWSBI),
$\eta=(\eta_x(t))_{x\in\Z, t\ge 0}$. Let $\widetilde \eta^{(+\varepsilon)}$ be the Poisson system of random walks
introduced in Section~\ref{Sect:Poissyst}. We then attempt to add extra particles (labeled as $\widehat \eta$ particles) to the Poisson system $\widetilde \eta^{(+\varepsilon)}$ at the origin with rate $\gamma$ provided that the origin is vacant under $\widetilde \eta^{(+\varepsilon)}$, and these attempted additions are coupled with those in the $\eta$ system. In particular, a particle added in the $\eta$ system can be coupled either to an $\widehat\eta$ particle added at the same time if the origin is vacant under $\widetilde \eta^{(+\varepsilon)}$, or to a particle in the $\widetilde\eta^{(+\varepsilon)}$ system if the origin is occupied under $\widetilde \eta^{(+\varepsilon)}$. This coupling constructs the $\eta$ particles as a subset of the $\widetilde \eta^{(+\varepsilon)}$ and $\widehat \eta$ particles, for which explicit calculations can be carried out.

\subsection{Coupling with the Poisson system}

We now formulate precisely the coupling between the Poisson system $\widetilde{\eta}^{(+\varepsilon)}$, the system of particles $\widehat \eta^{(+\varepsilon)}$ added during the times when $\widetilde{\eta}^{(+\varepsilon)}$ is vacant at the origin, and the true RWSBI system $\eta$.

Suppose that the Poisson system $\widetilde{\eta}^{(+\varepsilon)}$ has been constructed. Let $0<T_1<T_2< \cdots $ be the times of an independent rate $\gamma$ Poisson point process on $[0,\infty)$. At each time $T_i$, we add a particle at the origin to the $\widehat \eta^{(+\varepsilon)}$ system if the origin is vacant under $\widetilde \eta^{(+\varepsilon)}$. The successfully added particles then perform independent random walks. We now construct the $\eta$ system from $\widetilde \eta^{(+\varepsilon)}$ and $\widehat \eta^{(+\varepsilon)}$ as follows.

\begin{itemize}
\item  At time $T_1$, the origin is either occupied by a particle in the Poisson system $\widetilde{\eta}^{(+\varepsilon)}$, or a particle is added at the origin to the $\widehat{\eta}^{(+\varepsilon)}$ system. In either case, we add a particle to $\eta$ at the origin, which follows the same random walk as the particle (pick one if there is more than one) at the origin in the union of $\widetilde{\eta}^{(+\varepsilon)}$ and $\widehat{\eta}^{(+\varepsilon)}$.

\item  Assume that by time $T_{k}$ for some $k\geq 1$, particles have been added to $\eta$ in such a way that each particle in $\eta$ is coupled to a distinct particle in the union of $\widetilde{\eta}^{(+\varepsilon)}$ and $\widehat{\eta}^{(+\varepsilon)}$. We now attempt to add a particle at time $T_{k+1}$ to $\eta$ that preserves this coupling condition.

    \begin{itemize}
        \item If the origin is occupied at time $T_{k+1}$ under $\eta$, then no particle is added to $\eta$.

        \item If the origin is vacant at time $T_{k+1}$ under $\eta$, we note that it is either occupied under the Poisson system $\widetilde{\eta}^{(+\varepsilon)}$, or a particle is added at the origin to the $\widehat{\eta}^{(+\varepsilon)}$ system. In either case, the origin is occupied by particles in the union of $\widetilde{\eta}^{(+\varepsilon)}$ and $\widehat{\eta}^{(+\varepsilon)}$, and none of these particles could have been coupled with any particle in $\eta$. We then add a particle at the origin to $\eta$, which follows the same random walk as a corresponding particle in the union of $\widetilde{\eta}^{(+\varepsilon)}$ and $\widehat{\eta}^{(+\varepsilon)}$ at the origin.
    \end{itemize}
\end{itemize}

From the above inductive construction of $\eta$, it is clear that each particle in $\eta$ is coupled to a distinct particle in the union of $\widetilde{\eta}^{(+\varepsilon)}$ and $\widehat{\eta}^{(+\varepsilon)}$, and hence almost surely,
\begin{align}\label{domination}
\eta_x(t) & \leq \widetilde{\eta}^{(+\varepsilon)}_x(t) +
\widehat{\eta}^{(+\varepsilon)}_x(t) \qquad \mbox{for all } x \in \Z, \, t\geq 0 \\
\intertext{and in particular}
\label{dominationsum}
\sum_x \eta_x(t) & \leq \sum_x \widetilde{\eta}^{(+\varepsilon)}_x(t) + \sum_x \widehat{\eta}^{(+\varepsilon)}_x(t) \qquad \mbox{for all } t\geq 0.
\end{align}

\subsection{Proof of Theorem~\ref{thm:growth} (upper bound)}
By (\ref{dominationsum}), for any $\varepsilon>0$, we have
$$
\limsup_{t\to\infty} \frac{1}{R(t)} \sum_x \eta_x(t) \leq \limsup_{t\to\infty} \frac{1}{R(t)} \sum_x \widetilde{\eta}^{(+\varepsilon)}_x(t) + \limsup_{t\to\infty} \frac{1}{R(t)} \sum_x \widehat{\eta}^{(+\varepsilon)}_x(t),
$$
where the first term equals $1+\varepsilon$ by (\ref{eq:etaplus.totalno}). The second term equals $0$ because by construction, conditioned on $\widetilde{\eta}^{(+\varepsilon)}$, $\sum_x \widehat{\eta}^{(+\varepsilon)}_x(t)$ is a time-changed Poisson process with mean $\gamma\widetilde{V}^{(+\varepsilon)}_{0,t}$, and $\widetilde{V}^{(+\varepsilon)}_{0,t}/R(t)\to 0$ a.s.\ as $t\to\infty$ by (\ref{eq:etaplus.vacant}).
Therefore
$$
\limsup_{t\to\infty} \frac{1}{R(t)} \sum_x \eta_x(t) \leq 1+\varepsilon,
$$
which gives the desired upper bound if we let $\varepsilon\downarrow 0$.
\qed

\subsection{Proof of Theorem~\ref{thm:shape} (upper bound)} \label{subsec:shapeub}
By (\ref{domination}), for any $\varepsilon>0$ and any bounded non-negative continuous function $f\in C_{b,+}(\R)$, we have
\begin{align}
\frac{1}{\sigma\sqrt{t} \log t}
\sum_{x\in\Z} \eta_{x}(t) f\big(\tfrac{x}{\sigma\sqrt{t}}\big)
\le \frac{1}{\sigma\sqrt{t} \log t}
\sum_{x\in\Z} \widetilde{\eta}^{(+\varepsilon)}_x(t) f\big(\tfrac{x}{\sigma\sqrt{t}}\big)
+ || f ||_\infty \frac{\sum_x \widehat{\eta}^{(+\varepsilon)}_x(t)}{\sigma\sqrt{t} \log t}.
\end{align}
Since $R(t) \sim \sigma(\frac{2}{\pi})^{\frac{1}{2}} \sqrt{t} \log t$ by (\ref{rhod1totalgrowth}), the second term tends to $0$ as $t\to\infty$ as shown above in the proof of Theorem~\ref{thm:growth}, and hence almost surely,
\begin{align}
\label{eq:ub.profile}
\limsup_{t\to\infty} \frac{1}{\sigma\sqrt{t} \log t}
\sum_{x\in\Z} \eta_{x}(t) f\big(\tfrac{x}{\sigma\sqrt{t}}\big)
\leq \limsup_{t\to\infty} \frac{1}{\sigma\sqrt{t} \log t}
\sum_{x\in\Z} \widetilde{\eta}^{(+\varepsilon)}_x(t) f\big(\tfrac{x}{\sigma\sqrt{t}}\big).
\end{align}
Denote $\Xi_t:= \sum_{x\in\Z} \widetilde{\eta}^{(+\varepsilon)}_x(t) f\big(\tfrac{x}{\sigma\sqrt{t}}\big)$.
First we note that
$$
\frac{\E[\Xi_t]}{\sigma\sqrt{t} \log t}
= (1+\varepsilon)\sum_{x\in \Z} \frac{\rho_x(t)}{\sigma\sqrt{t}\log t} f\big(\tfrac{x}{\sigma\sqrt{t}}\big) \underset{t\to\infty}{\longrightarrow}
(1+\varepsilon) \int_\R f(y) \tilde \rho(y) dy =: M,
$$
where the convergence follows from Lemma~\ref{lem:rhoscaling} and a Riemann sum approximation of the integral. To show that $\Xi_t/(\sigma\sqrt{t}\log t)$
converges a.s.\ to the same limit $M$, we note that $\Xi_t$ is a weighted sum of independent Poisson random variables with mean $m_t:=\E[\Xi_t]=(M+o(1))\sigma\sqrt{t}\log t$, and each individual weight is uniformly bounded by $|| f ||_\infty$. By elementary large deviation estimates for Poisson random processes, for any $\delta>0$, we have
$$
\P(|\Xi_t-m_t|\geq \delta m_t) \leq C_1 e^{-C_2 m_t} \leq C_1 e^{-C_3 \sqrt{t}\log t},
$$
and hence by Borel-Cantelli, $\Xi_t/m_t\to 1$ a.s.\ along the time sequence $t_n=(\log n)^2$. To extend it to all $t\uparrow \infty$, by Borel-Cantelli, it suffices to show that for each $\delta>0$,
\begin{equation}\label{flucbd}
\sum_{n=1}^\infty \P\Big(\sup_{t\in [t_n, t_{n+1}]} |\Xi_t-\Xi_{t_n}| \geq \delta m_{t_n} \Big) <\infty.
\end{equation}
Note that $t_{n+1}-t_n\sim 2\log n/n$, and $\sup_{t\in [t_n, t_{n+1}]} |\Xi_t-\Xi_{t_n}|$ can be bounded in terms of the number of particles added to the $\widetilde \eta^{(+\varepsilon)}$ system during the time interval $[t_n, t_{n+1}]$ (which is Poisson distributed), plus the number of particles in $\widetilde \eta^{(+\varepsilon)}_{t_n}$ which have unusually large displacements (of order $\sqrt{t_{n}}$) during $[t_n, t_{n+1}]$ (note that these displacements are independent). Elementary large deviation estimates then give (\ref{flucbd}).

In conclusion, the RHS of (\ref{eq:ub.profile}) converges a.s.\ to $(1+\varepsilon) \int_\R f(y) \tilde{\rho}(y)\, dy$. Since $\varepsilon>0$ can be arbitrary, this implies the desired upper bound in Theorem~\ref{thm:shape}.
\qed

% %%%%%%%%%%%%%%%%%%%%%%%%%%% Upper bd %%%%%%%%%%%%%%%%%%%%%%%%%%

\section{Lower bounds in Theorems~\ref{thm:growth} and \ref{thm:shape}} \label{Sect:lb}

Here is the basic strategy for the lower bound on the system of random walks with self-blocking immigration (RWSBI),
$\eta=(\eta_x(t))_{x\in\Z, t\ge 0}$. Let $\widetilde \eta:=\widetilde \eta^{(-\varepsilon)}$ be a Poisson system of random walks with immigration rate $\beta^{(-\varepsilon)}$ as introduced in Section~\ref{Sect:Poissyst}. To get a lower bound on the $\eta$ system, we will construct an auxiliary system of $\widehat \eta$ particles, where particles are added at rate at most $\gamma$ and only when the origin is vacant under $\widehat \eta$, and $\widehat\eta$ particles may be killed from time to time. Such an $\widehat \eta$ system will be embedded as a subset of the $\eta$ system. To have explicit control on the rate at which particles are added in the $\widehat \eta$ system, which will lead to a lower bound on $\eta$, we couple $\widehat \eta$ with the Poisson system $\widetilde \eta$ in such a way that each particle added to $\widehat \eta$ is coupled with a particle in $\widetilde \eta$ (albeit starting at a different time), so that when we attempt to add a new particle to $\widehat \eta$, the origin being vacant under $\widetilde \eta$ ensures that it is also vacant under $\widehat \eta$. We can then bound from below the rate at which $\widehat\eta$ particles are added  in terms of the vacant time (at the origin) of the Poisson system $\widetilde \eta$, which can be estimated explicitly. This strategy will be made more precise in the following subsections.

\subsection{Coupling of one-dimensional random walks}
\label{subsect:vaccoup.1part}

We will need the following result, which shows that for two random walks $X$ and $Y$ starting respective at $x\in\Z$ and $0$ at time $0$ with $|x|\gg 1$,  there is a coupling between $X$ and $Y$ such that with high probability, the coupling is successful in the sense that $X$ and $Y$ coalesce and become a single walk before time $\tau_0:=\inf\{t\geq 0:X(t)=0\}$. Furthermore, whether the coupling is successful or not is independent of $(X(\tau_0+t))_{t\geq 0}$.

\begin{lem}\label{L:coupling}
For $n\in\N$, let $X_n$ and $Y_n$ be two rate $1$ continuous time random walks on $\Z$ with increment distribution $(a_x)_{x\in\Z}$ as specified in \eqref{kernela}, starting respectively at $x_n$ and $0$ at time $0$. Then there exists a coupling between $X_n$ and $Y_n$ with a coupling time $T_n$, such that: 
\begin{itemize}
\item[\rm (i)]  Either $T_n \leq \tau^{X_n}_0:=\inf\{t\geq 0: X_n(t)=0\}$ and $X_n(t)=Y_n(t)$ for all $t\geq T_n$, in which case we call the coupling successful; or $T_n=\infty$ and we call the coupling unsuccesful;
\item[\rm (ii)] The event $F_n:=\{T_n\leq \tau^{X_n}_0\}$ is measurable w.r.t.\ $Y_n$ and $(X_n(t))_{0\leq t\leq \tau^{X_n}_0}$, and on its complement $\{T_n=\infty\}$, $(X_n(\tau^{X_n}_0+t))_{t\geq 0}$ is independent of $Y_n$ and $(X_n(t))_{0\leq t\leq \tau^{X_n}_0}$.
\item[\rm (iii)] If $|x_n|\to\infty$ as $n\to\infty$, then $\P(F_n)\to 1$.
\end{itemize}
Furthermore, when $p$ is symmetric, the coupling can be chosen such that the joint dynamics of $(X_n(t), Y_n(t))_{t \ge 0}$ is Markovian.
\end{lem}

\begin{rmk} \rm When $X_n$ and $Y_n$ are simple symmetric random walks on $\Z$, there is a simple Markovian coupling such that the coupling is successful with probability $1$ for all $n \in \N\,$: 

If $X_n(0)$ is even, then we let $X_n$ and $Y_n$ jump simultaneously 
but in opposite directions until the first time that the two walks meet, 
and from this time on they perform identical jumps. This ensures that $X_n$ and $Y_n$ coalesce before $X_n$ hits $0$. If $X_n(0)$ is odd, then we wait for the first jump by either $X_n$ or $Y_n$, when the difference becomes even, and then couple as before. 
\end{rmk}

\noindent
{\bf Proof.} Without loss of generality, we may assume that $x_n\to\infty$. 

When $p$ is symmetric, we can couple $X_n$ and $Y_n$ such that they take opposite steps (simply putting $Y_n(t):=x_n-X_n(t)$) until they get close (i.e., either they meet or exchange order), and then run them as independent random walks until either they meet or $X_n$ hits $0$, whichever happens first. In the first case we set the meeting time to be $T_n$ and let the two walks move together afterwards; in the second case we just set $T_n=\infty$. 
\smallskip

In the general case, we can still couple  $X_n$ and $Y_n$ such that they take ``essentially'' opposite steps until they get close 
by a suitable coupling to Brownian motion, and then proceed as above. 
To implement this strategy, let
$$
\tau^{Y_n}_{1/2}:=\inf\{t\geq 0: Y_n(t)\ge x_n/2\} \quad \mbox{and} \quad \tau^{X_n}_{1/2}:=\inf\{t\geq 0: X_n(t)\le x_n/2\}.
$$
By Donsker's invariance principle, as $x_n\to \infty$, 
\begin{equation}
\label{eq:invprincXY}
\begin{aligned}
((x_n^{-1} Y_n(x_n^2t))_{t\geq 0}, x_n^{-2} \tau^{Y_n}_{1/2}) & \Asto{n} ((B_t)_{t \geq 0}, \tau_{1/2}) \\
((x_n^{-1} X_n(x_n^{2}t))_{t\geq 0}, x_n^{-2} \tau^{X_n}_{1/2}) & \Asto{n} ((1- B_t)_{t\geq 0}, \tau_{1/2}),
\end{aligned}
\end{equation}
where $(B_t)_{t\geq 0}$ is a Brownian motion with $\E[B_t^2]= \sigma^2 t$ and $\tau_{1/2}:= \inf\{t\geq 0: B_t \geq 1/2\}$. By Skorohod's representation theorem, we can couple $(X_n)_{n\geq 1}$ and $B$, and also $(Y_n)_{n\geq 1}$ and $B$, first possibly on different probability spaces, such that in both lines of \eqref{eq:invprincXY} the convergence holds almost surely. Then, using regular versions of the conditional distribution given $B$ on both probability spaces together with the same Brownian motion, we can construct copies of $(X_n)_{n\geq 1}$, $(Y_n)_{n\geq 1}$ and $B$ on the same probability space such that the convergence in both lines of \eqref{eq:invprincXY} holds simultaneously almost surely. We will use this coupling, which forces $X_n$ and $Y_n$ to take essentially opposite steps.

Since $\tau^{X_n}_{1/2}$ and $\tau^{Y_n}_{1/2}$ are stopping times, we may resample $(X_n(t))_{t\geq \tau^{X_n}_{1/2}}$ and $(Y_n(t))_{t\geq \tau^{Y_n}_{1/2}}$ 
(conditional on $X_n(\tau^{X_n}_{1/2})$, respectively on $Y_n(\tau^{Y_n}_{1/2}$)) independently of their past and of each other without changing the law of $X_n$, resp.\ $Y_n$. Assume this resampling from now on, and let
$$
\tau_{X_n, Y_n} := \inf\{t \geq \tau^{X_n}_{1/2} \vee \tau^{Y_n}_{1/2}: X_n(t)=Y_n(t)\}.
$$

On the event $\tau_{X_n, Y_n}\leq \tau^{X_n}_0=\inf\{t\geq 0: X_n(t)=0\}$, we set the coupling time $T_n= \tau_{X_n, Y_n}$
and resample $Y_n$ to be equal to $X_n$ from time $T_n$ onward. The coupling is then successful.

On the event $\tau^{X_n}_0 < \tau_{X_n, Y_n}$,  we set $T_n=\infty$ and the coupling is unsuccessful, and we do not make any further modification of $X_n$ and $Y_n$. 

With the above coupling, properties (i) and (ii) in Lemma~\ref{L:coupling} are clearly satisfied. To verify (iii), we need to show that under our coupling,
\begin{equation}\label{xymeet}
\P(\tau_{X_n, Y_n} \leq \tau^{X_n}_0) \to 1 \qquad \mbox{as}\quad |x_n|\to \infty. 
\end{equation}
Note that the above probability does not change if we assume $(X_n(t))_{t\geq \tau^{X_n}_{1/2}}$ and $(Y_n(t))_{t\geq \tau^{Y_n}_{1/2}}$ are coalescing random walks starting respectively at the space-time points $(X_n(\tau^{X_n}_{1/2}), \tau^{X_n}_{1/2})$ and $(Y_n(\tau^{Y_n}_{1/2}), \tau^{Y_n}_{1/2})$, where under our coupling, $(x_n^{-1} X_n(\tau^{X_n}_{1/2}), x_n^{-2} \tau^{X_n}_{1/2})$ and $(x_n^{-1} Y_n(\tau^{Y_n}_{1/2}), x_n^{-2} \tau^{Y_n}_{1/2})$ converge almost surely to the same space-time point $(1/2, \tau_{1/2})$. Therefore by the weak convergence of coalescing random walks to coalescing Brownian motions (proved in~\cite[Section 5]{NRS05} for discrete time random walks and is easily seen to hold also in continuous time), $(x_n^{-1} X_n(x_n^{2}t))_{t\geq x_n^{-2} \tau^{X_n}_{1/2}}$ and $(x_n^{-1} Y_n(x_n^{2}t))_{t\geq x_n^{-2} \tau^{Y_n}_{1/2}}$ converge to the same Brownian motion $W$ starting at $1/2$ at time $\tau_{1/2}$, and the rescaled time of coalescence, $x_n^{-2} \tau_{X_n,Y_n}$, converges to $\tau_{1/2}$. In particular, the probability that $(X_n(t))_{t\geq \tau^{X_n}_{1/2}}$ hits $0$ before meeting $(Y_n(t))_{t\geq \tau^{Y_n}_{1/2}}$ tends to $0$ as $n$ tends to infinity. This implies \eqref{xymeet} and the claim in (iii).
\qed

\subsection{Coupling with the Poisson system}\label{subsec:coupling}

We now formulate precisely the coupling between the true system $\eta$, the Poisson system $\widetilde \eta:=\widetilde \eta^{(-\varepsilon)}$ with immigration rate $\beta^{(-\varepsilon)}(t) = (1-\varepsilon) \gamma e^{-\rho_0(t)}$, and the auxiliary system $\widehat \eta$ as outlined at the start of this section. To simplify notation, in the remainder of the subsection, we will drop $(-\varepsilon)$ from the superscript and simply write $\widetilde \eta$ instead of $\widetilde \eta^{(-\varepsilon)}$.

Let $t_0=0 < t_1 < t_2 < t_3< \cdots$, and consider the time intervals $I_n=(t_{3n-3},t_{3n-2}]$, $I_n'=(t_{3n-2},t_{3n-1}]$, $I_n''=(t_{3n-1},t_{3n}]$. The precise values of the $t_n$'s will be determined later in \eqref{tns}, with $|I_n|=|I_n'|\ll |I_n''|$.
We will attempt to add exactly one $\widehat \eta$ particle in each time interval $I_n$, which will be coupled with the first $\widetilde \eta$ particle added during the time interval $I_n''$, with the coupling prescribed in Lemma~\ref{L:coupling}.

More precisely, let $(\widetilde N_t)_{t\geq 0}$ be a Poisson process with rate $\beta^{(-\varepsilon)}$, which determine the times when particles are added to $\widetilde \eta$, and let $(N_t)_{t\geq 0}$ be an independent Poisson process with rate $\gamma$, which determines the times when we might attempt to add particles to $\widehat \eta$ and $\eta$. Start with $\widetilde \eta_\cdot(0)=\eta_\cdot(0)=\widehat{\eta}_\cdot(0) \equiv 0$, and as an inductive hypothesis, assume that particles have been added to $\widetilde \eta$, $\eta$ and $\widehat \eta$ up to time $t_{3(n-1)}$ for some $n\geq 1$, such that the following properties hold:
\begin{itemize}
\item[\rm (a)] The paths of all added $\widetilde \eta$ particles have been sampled to time $\infty$, while the path of each added $\eta$ 
particle has been sampled till its first return to the origin after time $t_{3(n-1)}$, and the same for each $\widehat \eta$ particle unless it has been killed earlier;
\item[\rm (b)] Each $\widehat\eta$ particle is coupled to a distinct $\eta$ particle so that they follow the same path till the time of death of the $\widehat \eta$ particle. In particular, there are always more $\eta$ particles at the origin than $\widehat \eta$ particles;
\item[\rm (c)] Each $\widehat \eta$ particle added during the time interval $I_k$, for any $k\leq n-1$, is either killed at its first return to the origin after time $t_{3k}$, or it lives forever and is successfully coupled as in Lemma~\ref{L:coupling} to an $\widetilde \eta$ particle added during the subsequent time interval $I_k''$;
\item[\rm (d)] As a consequence of (c), at any time $t\geq  t_{3(n-1)}$, if one of the $\widehat \eta$ particles added before time $t_{3(n-1)}$ is at the origin, then so is one of the $\widetilde \eta$ particles added before time $t_{3(n-1)}$.
\end{itemize}
We now add particles to $\widetilde \eta$, $\eta$ and $\widehat \eta$ in the time intervals $(t_{3(n-1)}, t_{3n}]$ as follows.
\begin{itemize}
\item Add particles to $\widetilde \eta$ during the time interval $I_n\cup I_n'$ according to the Poisson process $\widetilde N$, with particle trajectories sampled to time $\infty$ according to independent random walks. Evolve existing $\eta$ and $\widehat \eta$ particles further till their first return to the origin after time $t_{3n}$. 
\item Let
$$
\widehat T_n:= \inf \{ t\in I_n: \widetilde \eta_0(t)=0, \Delta N_t=1\}, \qquad \widetilde T_n:= \inf \{ t\in I_n'': \Delta \widetilde N_t=1\},
$$
where $\inf \emptyset := \infty$. If $\widehat T_n=\infty$, then no $\widehat\eta$ particles are added during $(t_{3(n-1)}, t_{3n}]$, and  $\eta$ and $\widetilde \eta$ particles are added independently according to their own rules untill time $t_{3n}$, and their paths are sampled such that property (a) above continues to hold by $t_{3n}$;
\item If $\widehat T_n<\infty$, then we add an $\widehat\eta$ particle at the origin at time $\widehat T_n$. If the origin is occupied  in $\eta$ at time $\widehat T_n$, then we let the added $\widehat \eta$ particle follow the same path $X_n$ as one of the $\eta$ particles at the origin chosen at random. If the origin is vacant in $\eta$ at time $\widehat T_n$, then we also add an $\eta$ particle at time $\widehat T_n$ and let both particles follow the same random walk $X_n$, sampled independently of everything else till its first return to the origin after time $t_{3n}$. (Should the $\widehat \eta$ particle be killed later in the construction, we understand that the $\eta$ particle will be unaffected.)
\item Continue to add $\eta$ and $\widetilde \eta$ particles independently according to their own rules until time $t_{3n}\wedge \widetilde T_n$, and sample their paths so that property (a) continues to hold by $t_{3n}$. If $\widetilde T_n =\infty$, then kill the added $\widehat \eta$ particle at time $\tau_n:=\inf\{t \geq t_{3n}: X_n(t)=0\}$. 
\item If $\widetilde T_n<\infty$, then add an $\widetilde \eta$ particle at the origin at time $\widetilde T_n$ and sample its path $Y_n$ according to the conditional law of $Y_n$, conditioned on $(X_n(t))_{\widetilde T_n\leq t\leq \tau_n}$, so that $(X_n, Y_n)$ follows the law of the coupled random walks $(X_n, Y_n)$ in Lemma~\ref{L:coupling}. Denote
$$
E_n:=\{ \mbox{$\widehat T_n<\infty$, $\widetilde T_n<\infty$, and $X_n$ and $Y_n$ are coupled successfully as in Lemma~\ref{L:coupling}}\}.
$$
If the coupling is successful, then let the added $\widehat \eta$ particle live forever, otherwise kill it at time $\tau_n=\inf\{t \geq t_{3n}: X_n(t)=0\}$.  
\item Continue to add $\eta$ and $\widetilde \eta$ particles independently according to their own rules till time $t_{3n}$, with their trajectories sampled so that property (a) continues to hold by $t_{3n}$.
\end{itemize}
We note a subtle point in the above coupling, namely that we need to show that the $\widetilde \eta$ particles added at times $(\widetilde T_n)_{n\in\N}$ are indeed distributed as independent random walks. This is true because by construction, conditioned
on $X_n(\widetilde T_n)$ for $n\in\N$ with $\widehat T_n, \widetilde T_n<\infty$, the $\widehat \eta$ particle trajectories $(X_n(t))_{\widetilde T_n\leq t\leq \tau_n}$ are jointly independent, while the path $Y_n$ of each $\widetilde \eta$ particle coupled to $X_n$ depends only on $(X_n(t))_{\widetilde T_n\leq t\leq \tau_n}$ by Lemma \ref{L:coupling}.

It is clear that properties (a)--(d) above continue to hold after adding all particles up to time $t_{3n}$, and hence by induction, they hold for all time. In particular, by the coupling between $\eta$ and $\widehat \eta$, for all $n\in\N$, we have
\begin{align}
\label{eq:hatTTpr.lb}
\sum_{x\in\Z}\eta_x(t_{3n}) \geq \sum_{x\in\Z}\widehat \eta_x(t_{3n}) \geq \sum_{j=1}^n \mathbf{1}_{E_j} = \sum_{j=1}^n \mathbf{1}_{\{\widehat{T}_j < \infty, \widetilde{T}_j < \infty\}} \mathbf{1}_{E_j}.
\end{align}
To prove the lower bounds in Theorem~\ref{thm:growth} and \ref{thm:shape}, we will use the following choice of $(t_i)_{i\in\N}$:
\begin{equation}\label{tns}
t_{3n}= \varepsilon^2 \frac{n^2}{(\log (n \vee 3))^2} \quad \mbox{and} \quad  t_{3n+2}-t_{3n+1}=t_{3n+1}-t_{3n} = \varepsilon^2 (n+1)^{1-\varepsilon/2}, \qquad n\geq 0.
\end{equation}
The choice of $t_i$ is motivated by the fact that from
\eqref{rhod1totalgrowth}, the time until the $n$-th particle appears
in the true system should be of order $n^2/(\log n)^2$.
Note that \eqref{tns} implies
\begin{equation}\label{tns2}
\sum_{n: t_{3n}\leq t} 1 \underset{t\to\infty}{\sim} \big( \sqrt{t}\log t\big)/2\varepsilon.
\end{equation}

\begin{rmk} \rm When the random walk jump kernel $p(\cdot)$ is symmetric, we can use the
Markovian coupling of random walks guaranteed by Lemma~\ref{L:coupling} to construct 
the coupled $\widetilde \eta$, $\widehat \eta$ and $\eta$ particle systems jointly as a Markov process, with the use of labels to distinguish whether a particle is an 
$\widetilde \eta$, $\widehat \eta$ or $\eta$ particle, or it has multiple labels due to the coupling. 
 
  Spelling out the generator of such a system is straightforward, though lengthy, 
  so we do not make it explicit here. 
  Briefly, at a time $\widehat{T}_n<\infty$, if the origin is empty in 
  the $\eta$ system, we add a particle $X_n$ with label ``$\eta\&\widehat{\eta}$'';
  while if the origin is occupied in $\eta$ but empty in $\widehat{\eta}$, we 
  change the label of one of the $\eta$ particles to ``$\eta\&\widehat{\eta}$'' 
  (and call this the $X_n$ particle). The $X_n$ particle evolves as a free random walk until time $\widehat{T}_n \wedge t_{3n}$.
  If $\widehat{T}_n<t_{3n}$, then we add a particle $Y_n$ with label ``$\widetilde\eta$'', 
  and $X_n$ and $Y_n$ then evolve jointly as a Markov process according to the Markovian coupling from  Lemma~\ref{L:coupling} 
  until either they meet (at which time the particles merge and henceforth 
  evolve as a free random walk with label ``$\eta\&\widehat{\eta}\&\widetilde\eta$''), 
  or the $X_n$ particle visits the origin before meeting $Y_n$ (from this time 
  the $X_n$ particle evolves as a free random walk with label ``$\eta$'' and 
  the $Y_n$ particle evolves independently as a free random walk with label 
  ``$\widetilde\eta$''). If $\widehat{T}_n \ge t_{3n}$, then we change the label 
  of the $X_n$ particle to ``$\eta$'' at time $t_{3n}$. In between, all other 
  particles (with their labels) evolve independently, and additions of $\eta$, resp.\  
  $\widetilde{\eta}$ particles are executed according to their respective rules. 
\end{rmk}

\subsection{Proof of Theorem~\ref{thm:growth} (lower bound)}

First we note that the number of particles added to the Poisson system $\widetilde \eta$ during the time interval $I_n''$, which we denote by $\widetilde M_n$, is a Poisson random variable with mean
$$
\begin{aligned}
\E[ \widetilde M_n] = \int_{t_{3n-1}}^{t_{3n}} \beta^{(-\varepsilon)}(s)\, ds
= \ & \int_{t_{3n-1}}^{t_{3n}} \frac{\sigma(1-\varepsilon)}{\sqrt{2 \pi}}
\frac{\log s}{\sqrt{s}} \big(1 + o(1)\big) ds \\
\underset{n\to\infty}{\sim}\ &  \frac{2\sigma(1-\varepsilon)}{\sqrt{2 \pi}}
(\sqrt{t_{3n}} - \sqrt{t_{3n-1}}) \log t_{3n} \sim 4 \sigma\frac{\varepsilon(1-\varepsilon)}{\sqrt{2 \pi}},
\end{aligned}
$$
where we used the form of $\beta^{(-\varepsilon)}$ given in (\ref{betavareps}), the asymptotics for $\rho_0(t)$ given in (\ref{rhod1localgrowth}),
and the choice of $(t_i)_{i\in\N}$ given in (\ref{tns}). Therefore
$$
\P(\widetilde T_n<\infty) = \P(\widetilde M_n>0) = 1-e^{-\E[\widetilde M_n]} \underset{n\to\infty}{\longrightarrow} 1- e^{-4\sigma \frac{\varepsilon(1-\varepsilon)}{\sqrt{2 \pi}}}.
$$
Since $(\widetilde M_n)_{n\in\N}$ are independent, almost surely, we have
\begin{equation}\label{widetildeTnbd}
\sum_{j=1}^n \mathbf{1}_{\{\widetilde{T}_j < \infty\}} \underset{n\to\infty}{\sim} n\big(1- e^{-4 \sigma\frac{\varepsilon(1-\varepsilon)}{\sqrt{2 \pi}}}\big).
\end{equation}

Next we observe that on each time interval $I_n$, conditioned on the Poisson system $\widetilde \eta$,
$$
\P(\widehat T_n=\infty | \widetilde \eta) = e^{-\gamma \widetilde V_{t_{3n-3}, t_{3n-2}}},
$$
where by \eqref{eq:Evac.time.lb},
\begin{align*}
\E\big[ \widetilde{V}_{t_{3n-3}, t_{3n-2}}\big]
\ge c(t_{3n-2}-t_{3n-3})t_{3n-2}^{-(1-\varepsilon)/2} \ge
c\, \varepsilon^2 n^{1-\varepsilon/2} t_{3n}^{-(1-\varepsilon)/2}
= c\, \varepsilon^{1+\varepsilon} n^{\varepsilon/2} (\log n)^{1-\varepsilon}
\end{align*}
for some $c>0$. By the moment bound in Lemma~\ref{lem:lbPoiss.centr.k.mom} for $\big( \widetilde{V}_{t_{3n-3}, t_{3n-2}} - \E\big[\widetilde{V}_{t_{3n-3}, t_{3n-2}}\big] \big)^k$ for a sufficiently large even $k$ (note that the conditions are fulfilled), we can apply Markov's inequality and Borel-Cantelli to conclude that a.s., $\widetilde{V}_{t_{3n-3}, t_{3n-2}}/\E[\widetilde{V}_{t_{3n-3}, t_{3n-2}}]\to 1$, and hence $\{ \widetilde{V}_{t_{3n-3}, t_{3n-2}} > n^{\varepsilon/2}/2 \}$ occurs for all large enough $n$. Therefore a.s., $\sum_n \P(\widehat T_n=\infty | \widetilde \eta)<\infty$, and hence almost surely,
\begin{equation}\label{widehatTnbd}
\{\widehat T_j <\infty\} \quad \mbox{occurs for all } j \mbox{ sufficiently large}.
\end{equation}

Lastly we consider the events $E_j$ in (\ref{eq:hatTTpr.lb}). In our coupling construction of $\widetilde \eta$, $\eta$ and $\widehat \eta$, let ${\cal F}_n$ denote the $\sigma$-algebra generated by: the Poisson point process $\widetilde N$ up to time $\widetilde T_n \wedge t_{3n}$ and the trajectories of the $\widetilde\eta$ particles added before that time, as well as the Poisson point process $N$ up to time $\widehat T_n \wedge t_{3n-2}$ and the trajectories of the $\widehat \eta$ particles added before that time. Then $({\cal F}_n)_{n\in\N}$ defines a filtration, with $\{\widetilde T_n<\infty, \widehat T_n<\infty\}\in {\cal F}_n$, and $E_n \in {\cal F}_{n+1}$. Furthermore, since $|I_n'|\to\infty$, on the event $\{\widehat T_n<\infty, \widetilde T_n<\infty\}$, the path $X_n$ of the $\widehat \eta$ particle added at time $\widehat T_n$ satisfies $|X_n(\widetilde T_n)|\to\infty$ in probability as $n\to\infty$. Therefore by Lemma~\ref{L:coupling},
\begin{equation}\label{Enbd}
\begin{aligned}
& \big|\P(E_n | {\cal F}_n) - \mathbf{1}_{\{\widehat{T}_n < \infty, \widetilde{T}_n < \infty\}}\big| \\
\leq\ \, & \mathbf{1}_{\{\widehat{T}_n < \infty, \widetilde{T}_n < \infty\}} \big|\P(X_n \mbox{ and } Y_n \mbox{ are successfully coupled})-1\big|  \underset{n\to\infty}{\longrightarrow} 0.
\end{aligned}
\end{equation}
Note that (\ref{widetildeTnbd}) and (\ref{widehatTnbd}) imply
$$
\sum_{j=1}^n \mathbf{1}_{\{\widehat{T}_j < \infty, \widetilde{T}_j < \infty\}}  \underset{n\to\infty}{\sim} \sum_{j=1}^n \mathbf{1}_{\{\widetilde{T}_j < \infty\}}  \underset{n\to\infty}{\sim} n\big(1- e^{-4 \sigma\frac{\varepsilon(1-\varepsilon)}{\sqrt{2 \pi}}}\big) \longrightarrow  \infty \quad \mbox{a.s.},
$$
which together with (\ref{Enbd}) gives
\begin{equation}\label{Enbd2}
\sum_{j=1}^n \P(E_j | {\cal F}_j) \underset{n\to\infty}{\sim} \sum_{j=1}^n \mathbf{1}_{\{\widehat{T}_j < \infty, \widetilde{T}_j < \infty\}}  \underset{n\to\infty}{\sim} \sum_{j=1}^n \mathbf{1}_{\{\widetilde{T}_j < \infty\}}  \longrightarrow  \infty \quad \mbox{a.s.}
\end{equation}
On the other hand, by the second Borel-Cantelli Lemma~\cite[(4.11)]{D96},
$$
\frac{\sum_{j=1}^n \mathbf{1}_{E_j}}{\sum_{j=1}^n \P(E_j|{\cal F}_j)} \underset{n\to\infty}{\longrightarrow} 1  \quad \mbox{a.s.\  on}\quad  \Big\{\sum_{j=1}^\infty \P(E_j|{\cal F}_j)=\infty \Big\},
$$
which event is seen to have probability $1$ by (\ref{Enbd2}). Therefore, we also have
\begin{equation}\label{Enbd3}
\sum_{j=1}^n \mathbf{1}_{E_j} \underset{n\to\infty}{\sim} \sum_{j=1}^n \P(E_j | {\cal F}_j) \underset{n\to\infty}{\sim} \sum_{j=1}^n \mathbf{1}_{\{\widetilde{T}_j < \infty\}} \underset{n\to\infty}{\sim}  n\big(1- e^{-4\sigma \frac{\varepsilon(1-\varepsilon)}{\sqrt{2 \pi}}}\big)  \quad \mbox{a.s.}
\end{equation}
Since $t_{3j}= \varepsilon^2 \frac{j^2}{(\log j)^2}$, by (\ref{eq:hatTTpr.lb}) and \eqref{tns2}, this implies
\begin{equation}\label{etatlowerbd1}
\sum_x \eta_x(t) \geq \sum_{j: t_{3j}\leq t} \mathbf{1}_{E_j} \underset{t\to\infty}{\sim} \sum_{j=1}^{\sqrt{t}\log t/2\varepsilon} \mathbf{1}_{E_j} \underset{t\to\infty}{\sim} \frac{\sqrt{t}\log t}{2\varepsilon}\big(1- e^{-4\sigma \frac{\varepsilon(1-\varepsilon)}{\sqrt{2 \pi}}}\big) \quad \mbox{a.s.}
\end{equation}
Letting $\varepsilon\downarrow 0$ then gives the desired lower bound on $\sum_x \eta_x(t)$ in Theorem~\ref{thm:growth}.
\qed

\subsection{Proof of Theorem~\ref{thm:shape} (lower bound)}

The lower bound on the rate at which $\eta$ particles arrive readily leads to a lower bound on the spatial distribution of these particles at time $t$, since once an $\eta$ particle arrives, it evolves independently from all other particles.

First note that it suffices to verify Theorem~\ref{thm:shape} for $f\in C_{b,+}(\R)$ with a uniformly bounded derivative $f'$, since Theorem~\ref{thm:growth} on the convergence of the total mass of the measure $\frac{\sum_x \eta_x(t) \delta_x}{\sigma\sqrt{t}\log t}$ implies that it suffices to verify Theorem~\ref{thm:shape} for $f\in C_{b,+}(\R)$ with compact support, and any such $f$ can then be approximated in supremum norm by functions in $C_{b, +}(\R)$ with bounded derivatives.

Let us recall our construction of $\widehat \eta$ in Section~\ref{subsec:coupling}. For each $n\in\N$, an $\widehat\eta$ particle is added at time $\widehat T_n$ and then follows a random walk $X_n$ and lives forever, provided that $\widehat T_n<\infty$, $\widetilde T_n<\infty$, and $X_n$ can be successfully coupled to the random walk $Y_n$ that governs the motion of the $\widetilde \eta$ particle added at time $\widetilde T_n$. Then analogous to
(\ref{etatlowerbd1}), for any $f \in C_{b,+}(\R)$, almost surely
$$
\liminf_{t\to\infty}\frac{1}{\sigma{\sqrt t}\log t} \sum_x \eta_x(t) f\Big(\frac{x}{\sigma\sqrt t}\Big) \geq \liminf_{t\to\infty}\frac{1}{\sigma{\sqrt t}\log t} \sum_{n: t_{3n}\leq t} \mathbf{1}_{E_n} f\Big(\frac{X_n(t)}{\sigma\sqrt t}\Big).
$$
We can replace $\mathbf{1}_{E_n}$ above by $\mathbf{1}_{\{\widetilde T_n<\infty\}}$, since by (\ref{Enbd3}) and (\ref{tns2}), we have
$$
\begin{aligned}
& \frac{1}{\sigma{\sqrt t}\log t} \Big|\sum_{n: t_{3n}\leq t} \mathbf{1}_{E_n} f\Big(\frac{X_n(t)}{\sigma\sqrt t}\Big) -  \sum_{n: t_{3n}\leq t} \mathbf{1}_{\{\widetilde T_n<\infty\}} f\Big(\frac{X_n(t)}{\sigma\sqrt t}\Big)\Big| \\
=\ & \frac{1}{\sigma{\sqrt t}\log t} \sum_{n: t_{3n}\leq t} \!\! \big(\mathbf{1}_{\{\widetilde T_n<\infty\}} - \mathbf{1}_{E_n}\big) f\Big(\frac{X_n(t)}{\sigma\sqrt t}\Big)
\leq \frac{\Vert f\Vert_\infty}{\sigma{\sqrt t}\log t} \Big(\sum_{n: t_{3n}\leq t} \mathbf{1}_{\{\widetilde T_n<\infty\}} - \!\! \sum_{n: t_{3n}\leq t} \mathbf{1}_{E_n}\Big) \underset{t\to\infty}{\longrightarrow 0}.
\end{aligned}
$$
For any $\widehat \eta$ particle that gets killed, let us extend its path $X_n$ beyond its death by an independent random walk, and 
for $n$ with $\widehat T_n=\infty$, we let $X_n$ be an independent random walk starting from $t_{3n-2}$. We then have
\begin{equation}\label{fintegliminf}
\liminf_{t\to\infty}\frac{1}{\sigma{\sqrt t}\log t} \sum_x \eta_x(t) f\Big(\frac{x}{\sigma\sqrt t}\Big) \geq \liminf_{t\to\infty}\frac{1}{\sigma{\sqrt t}\log t} \sum_{n: t_{3n}\leq t} \mathbf{1}_{\{\widetilde T_n<\infty\}} f\Big(\frac{X_n(t)}{\sigma\sqrt t}\Big).
\end{equation}
Note that the space-time shifted random walks $(W_n(s))_{s\geq 0}:=(X_n(\widetilde T_n+s) -X_n(\widetilde T_n))_{s\geq 0}$, $n\in\N$, are i.i.d.\ and independent of the Poisson process $\widetilde N$ that determines the times when a particle is added to $\widetilde \eta$. We can then use $(W_n)_{n\in\N}$ and $\widetilde N$ to construct another Poisson system of random walks $\widetilde \xi$ with the same distribution as $\widetilde \eta$. More precisely, for each $n\in\N$ with $\widetilde T_n<\infty$, we add a $\widetilde \xi$ particle at the origin at time $\widetilde T_n$ which follows the trajectory $(W_n(s-\widetilde T_n))_{s\geq \widetilde T_n}=(X_n(s)-X_n(\widetilde T_n))_{s\geq \widetilde T_n}$. For all other times $t$ with $\Delta \widetilde N_t=1$, we add a $\widetilde \xi$ particle at the origin at time $t$, which follows an independent random walk trajectory.

We claim that a.s., the RHS of (\ref{fintegliminf}) does not change if we replace the trajectories of the $\widehat \eta$ particles therein by those of the $\widetilde \xi$ particles added at times $(\widetilde T_n)_{n\in\N}$. Indeed, the absolute difference arising from such a replacement (before taking $\liminf_{t\to\infty}$) is
\begin{align}
& \Big|\frac{1}{\sigma{\sqrt t}\log t} \sum_{n: t_{3n}\leq t} \mathbf{1}_{\{\widetilde T_n<\infty\}} \Big(f\Big(\frac{X_n(t)}{\sigma\sqrt t}\Big) - f\Big(\frac{X_n(t)-X_n(\widetilde T_n)}{\sigma\sqrt t}\Big)\Big)\Big| \notag\\
\leq \ & \frac{\Vert f'\Vert_\infty }{\sigma{\sqrt t}\log t} \sum_{n: t_{3n}\leq t}  \mathbf{1}_{\{\widetilde T_n<\infty\}} \min\Big\{\Vert f\Vert_\infty, \frac{|X_n(\widetilde T_n)|}{\sigma\sqrt t}\Big\} \notag\\
\leq\ & \frac{\Vert f'\Vert_\infty }{\sigma{\sqrt t}\log t} \sum_{n: t_{3n}\leq t} \min\Big\{\Vert f\Vert_\infty, \frac{\sup_{s\in [0, t_{3n}-t_{3n-3}]} |\widehat X_n(s)|}{\sigma\sqrt t}\Big\} , \label{fintegliminf2}
\end{align}
where $\widehat X_n$ is the random walk obtained from $X_n$ by shifting its starting time to $0$, which are i.i.d.\ and independent of $\widetilde N$. If we denote the minima in (\ref{fintegliminf2}) by $U_{n,t}$, then by Doob's $L^2$ maximal inequality for $\widehat X_n$,
$$
\E[U_{n,t}] \leq \frac{2}{\sigma\sqrt t} \E\big[\widehat X^2_{t_{3n}-t_{3n-3}}\big]^{\frac{1}{2}}= 2\sqrt{t_{3n}-t_{3n-3}} \leq \frac{C\sqrt n}{\sqrt{t}\log n}.
$$
Therefore using (\ref{tns2}),
\begin{equation}\label{Untbd1}
\frac{\Vert f'\Vert_\infty }{\sigma{\sqrt t}\log t} \sum_{n: t_{3n}\leq t} \E[U_{n,t}] \leq \frac{C' \Vert f'\Vert_\infty}{t\log t} (\sqrt{t}\log t)^{\frac{3}{2}} \underset{t\to\infty}{\longrightarrow 0}.
\end{equation}
Since $(U_{n,t})_{n\in\N}$ are independent random variables uniformly bounded by $\Vert f\Vert_\infty$, a standard fourth moment calculation applied to $\frac{\Vert f'\Vert_\infty }{\sigma{\sqrt t}\log t} \sum_{n: t_{3n}\leq t} \big(U_{n,t}-\E[U_{n,t}]\big)$, together with Markov inequality and Borel-Cantelli, show that this sequence converges a.s.\ to $0$ along the times $(t_{3N})_{N\in\N}$ (and hence also along $t\uparrow \infty$).
Together with (\ref{Untbd1}), this implies that the bound in (\ref{fintegliminf2}) converges a.s.\ to $0$ as $t\to\infty$, and hence we can replace the $\widehat \eta$ particle trajectories in the RHS of (\ref{fintegliminf}) by those of the $\widetilde \xi$ particles added at times $(\widetilde T_n)_{n\in\N}$.

We now make one more reduction, namely that including all particles in the $\tilde\xi$ system (not just those added at times $\widetilde T_n$) only introduces a small $\varepsilon$-dependent error. More precisely,
\begin{equation*}
\begin{aligned}
& \frac{1}{\sigma{\sqrt t}\log t} \Big|\sum_{n: t_{3n}\leq t} \mathbf{1}_{\{\widetilde T_n<\infty\}} f\Big(\frac{X_n(t)-X_n(\widetilde T_n)}{\sigma\sqrt t}\Big)
-\sum_x \widetilde \xi_x(t) f\Big(\frac{x}{\sigma\sqrt t}\Big)\Big| \\
\leq\ & \frac{\Vert f\Vert_\infty}{\sigma{\sqrt t}\log t} \Big|\sum_x \widetilde \xi_x(t) - \!\!\! \sum_{n: t_{3n}\leq t}\!\!\! 1_{\{\widetilde T_n<\infty\}}  \Big|
\underset{t\to\infty}{\longrightarrow} \Vert f\Vert_\infty \Big| \Big(\frac{2}{\pi}\Big)^{1/2}(1-\varepsilon)-\frac{1}{2\varepsilon \sigma}\big(1- e^{-4\sigma \frac{\varepsilon(1-\varepsilon)}{\sqrt{2 \pi}}}\big)  \Big| =: A_\varepsilon
\end{aligned}
\end{equation*}
by the a.s.\ asymptotics for $\sum_x \widetilde \xi_x(t)$ in Lemma~\ref{L:uppPoissonbd} and the a.s.\ asymptotics for $\sum_{n:t_{3n}\leq t} 1_{\{\widetilde T_n<\infty\}}$ given in (\ref{Enbd3}) and (\ref{etatlowerbd1}). Note that the limit $A_\varepsilon$ above tends to $0$ as $\varepsilon\downarrow 0$.

By the successive reductions we have made, we have thus shown that a.s.,
\begin{equation}\label{fintegliminf3}
\liminf_{t\to\infty}\frac{1}{\sigma{\sqrt t}\log t} \sum_x \eta_x(t) f\Big(\frac{x}{\sigma\sqrt t}\Big) \geq \liminf_{t\to\infty}\frac{1}{\sigma{\sqrt t}\log t} \sum_x \widetilde \xi_x(t) f\Big(\frac{x}{\sigma\sqrt t}\Big) - A_\varepsilon,
\end{equation}
where the above limit for the $\widetilde \xi$ system equals $(1-\varepsilon)\int_\R f(y) \tilde{\rho}(y)\, dy$ by the same argument as that in Section~\ref{subsec:shapeub} for the $\widetilde \eta^{(+\varepsilon)}$ system. Letting $\varepsilon\downarrow 0$ then completes the proof of the a.s.\ lower bound in Theorem~\ref{thm:shape}.
\qed

% %%%%%%%%%%%%%%%%%%%%%%%%%%% Appendix %%%%%%%%%%%%%%%%%%%%%%%%%%
\begin{appendix}

\section{Asymptotics of a semilinear lattice heat equation}
\label{Sect:rhoasymptotics}

This section is adapted from \cite[Sect.~3.6.1]{B03}
for ease of reference and the reader's convenience.
\smallskip

We consider the long-time behaviour of
the solution of the following inhomogeneous heat
equation on $\Z$ (which reduces to \eqref{eqnrho} upon
choosing $\alpha=1$):
\begin{eqnarray}
\label{eqnrhoapp}
\partial_t \rho_x(t) &=& L_{\mathrm{rw}} \rho_x(t)
        +\gamma \delta_0(x) \exp(-\alpha \rho_0(t)),
        \quad t \ge 0, \; x \in \Z, \\
\nonumber
\rho_x(0) & \equiv & 0,
\end{eqnarray}
where $\gamma, \alpha >0$ are parameters, and $L_{\rm rw}$ is the generator of a rate $1$ continuous time random walk $X$ on $\Z$, whose jump increments follow the probability kernel $(a_{-x})_{x\in\Z}$ with mean $0$ and variance $\sigma^2$, as specified in (\ref{kernela}).

\begin{rmk}
\label{remarkchoiceofgamma}
\rm
1.\ In integral form (sometimes called ``Duhamel's principle''),
\eqref{eqnrhoapp} reads
\begin{equation}
\label{eq:rhoDuhamel}
\rho_x(t) = \gamma \int_0^t p_{x}(t-s) \exp(-\alpha \rho_0(s)) \, ds,
\quad x \in \Z, \: t \ge 0
\end{equation}
where $p_x(t) = \P_0(X(t)=x)$ is the transition probability of
a continuous-time random walk with generator $L_{\mathrm{rw}}$.
\smallskip

\noindent 2.\ Let $\rho$ be the solution of \eqref{eqnrhoapp}.
Then $\vartheta_x(t) := \alpha \rho_x(t)$ solves
$\partial_t \vartheta_x(t) = L_{\mathrm{rw}}\vartheta_x(t)
        + \gamma' \delta_0(x) \exp(-\vartheta_0(t))$
with $\gamma':=\gamma \alpha$,
hence it suffices to consider the case $\alpha=1$.
\smallskip

\noindent 3.\ \eqref{eqnrhoapp} (and hence also \eqref{eqnrho})
has a unique solution:
Let $\rho^{(1)}$, $\rho^{(2)}$ be solutions, then
\begin{align*}
\frac{\partial}{\partial_t} \sum_x \big(\rho^{(2)}_x(t)-\rho^{(1)}_x(t) \big)^2
& = 2\sum_x \big(\rho^{(2)}_x(t)-\rho^{(1)}_x(t) \big)
L_{\mathrm{rw}}(\rho^{(2)}_\cdot(t)-\rho^{(1)}_\cdot(t) \big)_x \\
& {} \hspace{2em}
+ 2\gamma \big(\rho^{(2)}_0(t)-\rho^{(1)}_0(t) \big)
\big( e^{-\alpha \rho^{(2)}_0(t)} - e^{-\alpha \rho^{(1)}_0(t)} \big) \; \leq 0,
\end{align*}
noting that $\sum_x f_x (L_{\mathrm{rw}} f)_x \leq 0$ for any $f \in \ell^2(\Z)$,
and $(a-b) \big( e^{-\alpha a} - e^{-\alpha b} \big) \leq 0$ for any $a, b \in \R$.
Hence $\rho^{(1)} \equiv \rho^{(2)}$.
\end{rmk}

\begin{lem}
\label{lemmarhoasympt} 
Let $\rho$ be the solution of \eqref{eqnrhoapp}. Then $\rho_0(t)$ is increasing in $t$, and
as $t\to \infty$,
\begin{eqnarray}
\label{rhod1localgrowth-app}
\rho_0(t) & = & \frac{1}{\alpha}
        \left\{\frac{1}{2}\log t - \log\log t + \log\Big(\sqrt{2\pi}
        \, \gamma \alpha/\sigma \Big)
        \right\} + o(1), \hspace{0em} \\
\label{rhod1totalgrowth-app}
\sum_x \rho_x(t) & = & \gamma \int_0^t e^{-\alpha \rho_0(s)} ds \sim
        \frac{\sigma}{\alpha} \left(\frac{2}{\pi}\right)^{1/2} \sqrt{t} \log t.
\end{eqnarray}
\end{lem}
\smallskip

\noindent
{\bf Proof.}
Assume w.l.o.g.\ $\alpha = 1$, cf.\ Remark~\ref{remarkchoiceofgamma}.
We see from \eqref{eq:rhoDuhamel} for $x=0$ that $\rho_0(t)$ is the solution of
the functional equation
\begin{equation}
\label{integralequationrho0}
f(t) = \int_0^t \gamma p_0(t-s) \exp(-f(s)) \, ds, \quad t \ge 0.
\end{equation}

Let us call a function $\bar{\varphi}: \Z^d \times \R_+ \to \R_+$
with $\bar{\varphi}_\cdot(0) \equiv 0$
a {\em strict supersolution to \eqref{eqnrhoapp}\/}
if it solves
\begin{equation}
\begin{array}{l}
\displaystyle
\partial_t \bar{\varphi}_x(t) = L_{\mathrm{rw}} \bar{\varphi}_x(t)+
        \gamma r^{\bar{\varphi}}(t) \delta_0(x) ,
        \quad t \ge 0, x \in \Z \\[1.5ex]
\mbox{with an  $r^{\bar{\varphi}}(t) > \exp(-\bar{\varphi}_0(t))$.}
\end{array}
\end{equation}
Then we see that $\bar{\varphi}_0(t) \ge \rho_0(t)$ for all $t \ge 0$:
Indeed,
$\psi_x(t) := \bar{\varphi}_x(t) - \rho_x(t)$ solves
\[
 \partial_t \psi_x(t) =  L_{\mathrm{rw}} \psi_x(t)
                        +\gamma % \underbrace{%
        (r^{\bar{\varphi}}(t)-e^{-\rho_0(t)})\delta_0(x) %}_{=:r^{\psi}(t)}\delta_0(x)
\]
and $\psi_0(t) > 0$ for small $t$.
Assume that $t_0 := \inf\{ t : \psi_0(t) < 0 \} < \infty$.
Then we would have $\psi_0(t_0) = 0$ by continuity,
but also $\psi_x(t_0) \ge 0$ for all $x$. To see this
observe that $\psi_x(t)$, $x \neq 0$ has a representation
($\psi$ solves the heat equation away from $0$, consider
$\psi_0(t)$ as exogenous input)
\[
\psi_x(t) = \int_0^t \psi_0(t-s) \P_x(T_0 \in ds)
+ \E_x [\psi_{X(t)}(0); T_0 > t]
\quad \Big(= \E_x \big[ \psi_{X(t \wedge T_0)}(t-(t \wedge T_0))\big] \,\Big)
\]
where $T_0 := \inf \{ s : X_s=0\}$ (see Lemma \ref{lemmaheateqnrepr}).
Hence $\psi_x(t_0) \ge 0$ for all $x$ because $\psi(0) \equiv 0$
and $\psi_0(s) \ge 0$ for $0 \le s \le t_0$ by definition.
Consequently $L_{\mathrm{rw}} \psi_0(t_0) \ge 0$
and we conclude that $\gamma^{-1} \partial_t \psi_0(t_0) \ge %% r^{\psi}(t_0) :=
r^{\bar{\varphi}}(t_0)-e^{-\rho_0(t_0)} >
\exp(-\bar{\varphi}_0(t_0))-\exp(-\rho_0(t_0)) = 0$
in contradiction to the definition of $t_0$.
\smallskip

We can construct a supersolution to (\ref{eqnrhoapp})
from a strict subsolution to (\ref{integralequationrho0}):
Assume $\underline{f} : [0,\infty) \to \R$ satisfies
\begin{equation}
\label{ineqstrictsubsol0}
 \underline{f}(t) < \int_0^t \gamma p_0(t-s)
                                         \exp(-\underline{f}(s)) \, ds
        \quad \mbox{for $t\ge 0$}.
\end{equation}
Then
\begin{equation}\label{supphi}
\bar{\varphi}_x(t):=
\int_0^t \gamma p_x(t-s) \exp(-\underline{f}(s)) \, ds
\end{equation}
solves
\[
\partial_t \bar{\varphi}_x(t) = L_{\rm rw} \bar{\varphi}_x(t)
        + \gamma \exp(-\underline{f}(t)) \delta_0(x)
\]
and in particular $\bar{\varphi}_0(t) > \underline{f}(t)$,
hence $\exp(-\underline{f}(t)) > \exp(-\bar{\varphi}_0(t))$.

Similarly, if $\underline{\varphi}$ is a strict subsolution
we have $\underline{\varphi}_0(t) \le \rho_0(t)$ for all $t \ge 0$
and such a $\underline{\varphi}$ can be constructed
analogously from
a supersolution $\bar{f}$ to (\ref{integralequationrho0}).
\smallskip

Observe that the solution $\rho$ of (\ref{eqnrhoapp})
has the property that $\rho_0(t)$ is an
increasing function: Obviously $\partial_t \rho_0(t) > 0$
for $t$ small.
Assume that $t_0 := \inf \{ t : \partial_t \rho_0(t) < 0 \} < \infty$.
Then by continuity $\partial_t \rho_0(t_0) = 0$.
We have for $x \in \Z \setminus \{0\}$
by the representation given in Lemma \ref{lemmaheateqnrepr}
\begin{eqnarray*}
\partial_t \rho_x(t_0) & = & \lim_h \frac{1}{h} \left[
  \int_0^{t_0} \rho_0(t_0-s) \P_x(T_0 \in ds)
  - \int_0^{t_0-h} \rho_0(t_0-h-s) \P_x(T_0 \in ds)
\right] \\
& = & \lim_h  \int_0^{t_0-h}
\frac{1}{h} (\rho_0(t_0-s) - \rho_0(t_0-h-s)) \P_x(T_0 \in ds) \\
&& {} + \lim_h\frac{1}{h} \int_{t_0-h}^{t_0} \rho_0(t_0-s) \P_x(T_0 \in ds) \\
& \geq & \int_0^{t_0} \partial_t \rho_0(t-s) \P_x(T_0 \in ds)
+ \underbrace{\rho_0(0)}_{=0} \frac{\P_x(T_0 \in dt)}{dt}|_{t=t_0}
> 0,
\end{eqnarray*}
because $\partial_t \rho_0(t) > 0$ in $[0,t_0)$ and $\mbox{supp}({\cal L}_x(T_0)) = \R_+$,
and we applied Fatou's Lemma in the first inequality. Thus
\[
\partial^2_t \rho_0(t_0) =
\sum_{x} a_x (\partial_t \rho_x(t_0) -\partial_t \rho_0(t_0))
\; - \partial_t \rho_0(t_0) \gamma \exp(-\rho_0(t_0)) > 0,
\]
contradicting the definition of $t_0$.

\begin{lem}
\label{lemmaapproxsolutions}
Assume $\alpha=1$. \\ 
i) For $C < \log \left(\sqrt{2 \pi}\, \gamma/\sigma \right)$
there exists a $K>0$ such that
\[
\underline{f}(t) := \left\{ \begin{array}{cl}
        \frac{1}{2} \log t - \log\log t + C & \mbox{if} \; t \ge K, \\
        -1 & \mbox{if} \; 0 \le t < K
        \end{array} \right.
\]
is a strict subsolution for (\ref{integralequationrho0}).

\noindent
ii) For $C > \log \left( \sqrt{2 \pi}\, \gamma/\sigma \right)$
there exist $K, K' >0$ such that
\[
\bar{f}(t) := \left\{ \begin{array}{cl}
        \frac{1}{2} \log t - \log\log t + C & \mbox{if} \; t \ge K, \\
        K' & \mbox{if} \; 0 \le t < K
        \end{array} \right.
\]
is a strict supersolution for (\ref{integralequationrho0}).
\end{lem}
\noindent
{\bf Proof.}
This is a straightforward computation using the local central limit theorem,
$p_0(t) \sim (2\pi \sigma^2 t)^{-1/2}$.
Here are some details:\\[0.5ex]
{\it i)}.
Let $e^{-C} = (1+3\varepsilon)/\big( \sqrt{2 \pi}\, \gamma/\sigma \big)$
with $\varepsilon > 0$ small. We have $p_0(t) \ge
(1-\varepsilon)/\sqrt{2\pi \sigma^2 t}$ for $t \ge t_0(\varepsilon)$.
For $\underline{f}$ as in $i)$ and any $t \geq K \vee t_0(\varepsilon)$
we estimate
\begin{align*}
\int_0^{t} p_0(t-s) \gamma e^{-\underline{f}(s)} \, ds
& \ge \int_K^{t-K} p_0(t-s) \, e^{-C}\gamma \frac{\log s}{\sqrt{s}} \, ds \\
& \ge \frac{1+\varepsilon}{2\pi} \int_K^{t-K} \frac{\log s}{\sqrt{s (t-s)}} \, ds
= \frac{1+\varepsilon}{2\pi} \int_{K/t}^{1-K/t}
        \frac{\log t + \log u}{\sqrt{u(1-u)}} du \\
& \ge \frac{1+\varepsilon}{2\pi} \left\{
\log t \Big[ \int_0^1 \frac{1}{\sqrt{u(1-u)}} du
        - 4 \sqrt{2K/t} \Big] +
        \int_0^1 \frac{\log u}{\sqrt{u(1-u)}} du \right\}.
\end{align*}
Observing that $\int_0^1 (u(1-u))^{-1/2} \, du = \pi$ and
$\int_0^1 (u(1-u))^{-1/2} \log u \, du \in (-\infty,0)$ we see that
there exists $n \,(=n(\varepsilon)) \ge 1$ such that for all $K \ge 1$
\[
\int_0^t \gamma p_0(t-s) e^{-\underline{f}(s)}\, ds \ge
\frac{1+\varepsilon/2}{2} \log t > \underline{f}(t)
\quad \mbox{whenever $t \ge nK$}.
\]
On the other hand for $ t_0(\varepsilon) \le K < t < nK$ we have
\begin{eqnarray*}
\int_0^{t} \gamma p_0(s) \exp(-\underline{f}(t-s)) \, ds & \ge &
        \frac{\gamma (1-\varepsilon)}{\sqrt{2\pi \sigma^2}}
\int_{t-K}^t \frac{1}{\sqrt{u}} du
\ge \frac{\gamma (1-\varepsilon)}{\sqrt{2\pi \sigma^2}}
\int_{(n-1)K}^{nK} u^{-1/2} du \\
& = & \frac{\gamma (1-\varepsilon) \sqrt{2}}{\sigma \sqrt{\pi}}
\left(\sqrt{nK}-\sqrt{(n-1)K}\right)
\ge \frac{\gamma}{2\sigma} \sqrt{\frac{K}{2\pi n}}
\end{eqnarray*}
and $\underline{f}(t) \le \log (nK)$. So we just have to chose
$K \ge 1$ so big that $\frac{\gamma}{2\sigma}  \sqrt{\frac{K}{2\pi n}} >  \log (nK)$.
\smallskip

{\it ii)} can be treated similarly.
\qed

\medskip

\noindent
{\bf Proof of Lemma \ref{lemmarhoasympt}, continued.}\/
Constructing $\bar{\varphi}$ and $\underline{\varphi}$ as in \eqref{supphi}
from the functions $\underline{f}$ and $\bar f$ given in Lemma \ref{lemmaapproxsolutions},
with $\underline{\varphi} \le \rho \le \bar{\varphi}$, we see easily that
\begin{equation}
\label{roughasymptoticrho0}
 \rho_0(t) \sim \frac{1}{2} \log t \quad \mbox{as $t\to\infty$}.
\end{equation}
But we need a finer result, namely
$\rho_0(t) = \frac{1}{2} \log t - \log\log t + \log (\sqrt{2 \pi} \gamma/\sigma)+o(1)$.
We use Laplace transforms to strengthen the asymptotics
\eqref{roughasymptoticrho0}:

Denoting $\xi(t):=\gamma \exp(-\rho_0(t))$ we can write
(\ref{integralequationrho0}) as $\rho_0 = p_0 \ast \xi$,
after taking Laplace transforms this reads
\begin{equation}
\label{eqnrhoplaplacetransform}
\widehat{\rho_0}(\lambda) =
        \widehat{p_0}(\lambda) \widehat{\xi}(\lambda),
\quad \lambda > 0.
\end{equation}
We have $\widehat{p_0}(\lambda) \sim (2\sigma^2 \lambda)^{-1/2}$
as $\lambda \downarrow 0$.
From (\ref{roughasymptoticrho0}) and a Tauberian theorem
(see e.g.\ \cite{FeII}, Chap.~XIII.5, Thm.~4) %%, p.~423)
we conclude that
$\widehat{\rho_0}(\lambda) \sim \frac{1}{2\lambda} \log(1/\lambda)$,
hence $\widehat{\xi}(\lambda) \sim \sigma (2\lambda)^{-1/2} \log(1/\lambda)$
for $\lambda \downarrow 0$.
Invoking the Tauberian theorem in the other direction we get
\[
 \gamma \exp(-\rho_0(t))
        = \left( \sigma (2\pi)^{-1/2} t^{-1/2} \log t \right) (1+o(1)).
\]
(\ref{rhod1localgrowth-app}) follows by taking logarithms.
Observe that the use of the Tauberian theorem is justified because $\rho_0$, and hence also $\xi$, are monotone
functions. Finally observe that
$\int (\log s)/\sqrt{s} \,ds = 2 \sqrt{s} \log s - 4\sqrt{s}$
to obtain (\ref{rhod1totalgrowth-app}).
\hfill $\Box$

\begin{lem}
\label{lemmaheateqnrepr}
Let $\psi_{\cdot}(0) : \Z \to \R$ and $\psi_0(\cdot) : \R_+ \to \R$
be given real-valued continuous functions and define $\psi$
on $\Z \times \R_+$ as the solution of the heat equation
corresponding to $L_{\mathrm{rw}}$ away from $0$ with given boundary
behaviour, i.e., $\psi$ solves
\[
\partial_t \psi_x(t) =  L_{\mathrm{rw}} \psi_x(t), \quad
        x \in \Z \setminus \{0\}, t \geq 0.
\]
Then $\psi$ has the stochastic representation
\[
\psi_x(t) = \E_x \big[ \psi_{X(t \wedge T_0)}(t-(t \wedge T_0)) \big]
\]
where $(X(t))_{t\ge0}$ is a continuous-time random walk on $\Z$
with generator $L_{\mathrm{rw}}$ and
$T_0 := \inf \{ s > 0 : X(s) = 0 \}$ the hitting time of the
origin.
\end{lem}

\begin{lem} \label{lem:rhoscaling}
Let $\rho$ be the solution of \eqref{eqnrhoapp}. Then uniformly in $x\in \R$, we have the following convergence:
\begin{equation}
\label{eq:rhoscaling1}
\frac1{\log t}\rho_{[\sigma \sqrt{t} x]}(t) \mathop{\longrightarrow}_{t\to\infty}
\frac1{2\pi} \int_0^1 \frac1{\sqrt{s(1-s)}} e^{-x^2/(2s)}\, ds
=: \tilde{\rho}(x) = 1 - \Phi(|x|),
\end{equation}
where $\Phi(a):= \frac1{\sqrt{2\pi}} \int_{-\infty}^a e^{-z^2/2}\, dz$.
\end{lem}
\noindent
{\bf Proof.} By \eqref{rhod1localgrowth-app}, we have
$\exp(- \alpha \rho_0(t)) = \big(\sigma (\log t)/\gamma\sqrt{2\pi t}\big)(1+o(1))$,
and by the local CLT, $p_{[\sigma \sqrt{t} x]}(t) = (2\pi \sigma^2 t)^{-1/2} (e^{-x^2/2}+o(1))$ uniformly in $x\in \R$.
Thus, by (\ref{eq:rhoDuhamel}),
\begin{align*}
\frac1{\log t}\rho_{[\sigma\sqrt{t} x]}(t) & =
\frac{\gamma}{\log t}  \int_0^t \frac{\log s}{\gamma\sqrt{2\pi s}}
\frac1{\sqrt{2\pi (t-s)}} e^{-t x^2/2(t-s)}\, ds + o(1) \\
& = \frac{1}{2\pi} \int_0^t \frac{\log t + \log(s/t)}{\log t}
\big( (s/t)(1-s/t) \big)^{-1/2} e^{-x^2/2(1-s/t)} \, \frac{ds}{t} + o(1) \\
& = \frac{1}{2\pi} \int_0^1 \frac1{\sqrt{(1-u)u}} e^{-x^2/2u} \, du + o(1),
\end{align*}
where we substituted $u=1-s/t$ in the last line.

To prove the identity in \eqref{eq:rhoscaling1}, substituting $1/s = z$ and then $z-1=y$, we find
\begin{align*}
\frac{1}{2\pi} & \int_0^1 \frac{1}{\sqrt{s(1-s)}} e^{\frac{-x^2}{2s}} \, ds
= \frac{1}{2\pi} \int_1^\infty \left( z^{-1}(1-z^{-1})\right)^{-1/2}
        \exp\Big(\frac{-x^2}{2} z \Big) \frac{dz}{z^2}  \\
& = \frac{1}{2\pi} \int_1^\infty \frac{1}{z\sqrt{z-1}}
        \exp\Big(\frac{-x^2}{2} z\Big) dz
= \frac{1}{2\pi} e^{-x^2/2} \int_0^\infty \exp\Big(\frac{-x^2}{2} y\Big)
                \frac{1}{\sqrt{y}(y+1)} \, dy \\
& = \frac{1}{2\pi} e^{-x^2/2} \hat{f}(x^2/2),
\end{align*}
where $\hat{f}$ is the Laplace transform of $f(t)=1/((y+1)\sqrt{y})$.
A table of Laplace transforms (e.g.~\cite[29.3.114]{AS64}) shows that
$\hat{f}(z)=2 \pi e^z (1-\Phi(\sqrt{2 z}))$.
\qed

% %%%%%%%%%%%%%%%% Poisson vacant correlations %%%%%%%%%%%

\section{Correlation functions for Poisson vacant events}
\label{Sect:Poiscor}
In this section, we compute the correlation function for the events that a Poisson point process is vacant on each of $k$ given sets.
This is used to prove Lemma~\ref{lem:lbPoiss.centr.k.mom} on the centered moments of the origin's vacant time for a Poisson system of random walks.

\begin{lem}
\label{lem:Poisvac}
Let $(S, \cal B)$ be a measurable space, $\xi$ a Poisson point process on $S$ with intensity measure $\nu$.
Then for $k \in \N$, $E_1,\ E_2,\dots, E_k \in \cal S$ with $\nu(E_1), \dots, \nu(E_k) < \infty$, and $M \in \N\cup\{0\}$,
\begin{align}
\E\bigg[ \prod_{i=1}^k & \Big( \mathbf{1}(\xi(E_i)=0) - \P(\xi(E_i)=0)
\Big) \bigg] \notag \\
\label{eq:Poisvac.full}
& =
e^{-\sum_{i=1}^k \nu(E_i)} \sum_{n=1}^\infty \frac{1}{n!}
\sum_{\begin{array}{c} \scriptstyle I_1, \dots, I_n \subset \{1,\dots,k\} \\
\scriptstyle |I_1|, \dots, |I_n| \geq 2\\
\scriptstyle I_1 \cup \cdots \cup I_n = \{1,\dots,k\} \end{array}}
(-1)^{\sum_{j=1}^n |I_j|} \prod_{j=1}^n \nu\Big( \mathop{\cap}_{\ell \in I_j} E_\ell\Big)
\\
\label{eq:Poisvac.orderM}
& =
e^{-\sum_{i=1}^k \nu(E_i)} \sum_{n=1}^M \frac{1}{n!}
\sum_{\begin{array}{c} \scriptstyle I_1, \dots, I_n \subset \{1,\dots,k\} \\
\scriptstyle |I_1|, \dots, |I_n| \geq 2\\
\scriptstyle I_1 \cup \cdots \cup I_n = \{1,\dots,k\} \end{array}}
(-1)^{\sum_{j=1}^n |I_j|} \prod_{j=1}^n \nu\Big( \mathop{\cap}_{\ell \in I_j} E_\ell\Big)
\notag \\
& \qquad \qquad + e^{-\sum_{i=1}^k \nu(E_i)}
R_{M+1}\big( 2^k \max_{1 \leq i < j \leq k} \nu\big( E_i \cap E_j \big)\big),
\end{align}
where $| R_{M+1}(x) | \leq 2^k \frac{|x|^{M+1}}{(M+1)!} e^{|x|}$.
\end{lem}
\begin{rmk}{\rm Lemma~\ref{lem:Poisvac} allows us to control the $k$-point correlation function quantitatively in terms of
$\nu(E_i\cap E_j)$, $1\leq i<j\leq k$. This result should be well known, but we sketch the proof below for completeness and lack of a precise reference.
}
\end{rmk}
\noindent
{\bf Proof.}
Since $\P\big(\xi(B)=0\big)=e^{-\nu(B)}$ for any set $B\in \cal B$, we have
\begin{align}
&\E\bigg[ \prod_{i=1}^k \Big( \mathbf{1}(\xi(E_i)=0) - \P(\xi(E_i)=0)
\Big) \bigg] \notag \\
= \ & \sum_{I' \subset \{1,\dots,k\}} (-1)^{k-|I'|}
\P\Big( \xi\big( \mathop{\cup}_{\ell \in I'} E_\ell\big) = 0 \Big)
\prod_{j \not\in I'} \P\big(\xi(E_j)=0\big) \notag \\
= \ & \sum_{I' \subset \{1,\dots,k\}} (-1)^{k-|I'|}
\exp\Big[ -\nu\big( \mathop{\cup}_{\ell \in I'} E_\ell\big) -
\sum_{j \not\in I'} \nu(E_j) \Big]  \notag \\
= \ & e^{-\sum_{i=1}^k \nu(E_i)} \sum_{I' \subset \{1,\dots,k\}} (-1)^{k-|I'|}
\exp\Big[ -\nu\big( \mathop{\cup}_{\ell \in I'} E_\ell\big) +
\sum_{j \in I'} \nu(E_j) \Big] \notag \\
= \ & e^{-\sum_{i=1}^k \nu(E_i)} \sum_{I' \subset \{1,\dots,k\}} (-1)^{k-|I'|}
\exp\Big[ \sum_{I \subset I', |I| \geq 2} (-1)^{|I|}
\nu\big( \mathop{\cap}_{\ell \in I} E_\ell\big) \Big] \label{corrbd}
\end{align}
where we used the inclusion-exclusion principle in the last line.

Note that when we Taylor expand the rightmost exponential in (\ref{corrbd}), the zeroth order term is
$\sum_{I' \subset \{1,\dots,k\}} (-1)^{k-|I'|}=0$. For a fixed $I' \subset \{1,\dots,k\}$ and $n\in\N$, the $n$-th
order term of the Taylor expansion for the exponential is
\begin{align*}
& \sum_{I' \subset \{1,\dots,k\}} (-1)^{k-|I'|}
\frac{1}{n!} \sum_{\begin{array}{c} \scriptstyle I_1, \dots, I_n \subset I' \\
\scriptstyle |I_1|, \dots, |I_n| \geq 2\end{array}}
(-1)^{|I_1|+\cdots+|I_n|}
\prod_{j=1}^n \nu\big( \mathop{\cap}_{\ell \in I_j} E_\ell\big) \\
=\  &
\frac{1}{n!} \!\!\!\!\!\!\!\!\!\!\! \sum_{\begin{array}{c} \scriptstyle I_1, \dots, I_n \subset \{1,\dots,k\} \\
\scriptstyle |I_1|, \dots, |I_n| \geq 2\end{array}} \hspace{-2.5em}
(-1)^{\sum_{j=1}^k |I_j|}
\prod_{j=1}^n \nu\big( \mathop{\cap}_{\ell \in I_j} E_\ell\big)
\hspace{-1.5em}
\sum_{\begin{array}{c} \scriptstyle  I' \subset \{1,\dots,k\} \\
\scriptstyle I' \supset I_1 \cup \cdots \cup I_n \end{array}}
\hspace{-2em} (-1)^{k-|I'|}
= \frac{1}{n!} \!\!\!\!\!\!\!\!\!\!\!
\sum_{\begin{array}{c} \scriptstyle I_1, \dots, I_n \subset \{1,\dots,k\} \\
\scriptstyle |I_1|, \dots, |I_n| \geq 2\\
\scriptstyle I_1 \cup \cdots \cup I_n = \{1,\dots,k\} \end{array}}
\hspace{-2.5em}
(-1)^{\sum_{j=1}^n |I_j|} \prod_{j=1}^n \nu\Big( \mathop{\cap}_{\ell \in I_j} E_\ell\Big)
\end{align*}
since whenever $I_1 \cup \cdots \cup I_n \neq \{1,\dots,k\}$,
the summation over $I'$ gives $0$. This proves \eqref{eq:Poisvac.full}.

To check \eqref{eq:Poisvac.orderM}, let
$\phi(I') := \sum_{I \subset I', |I| \geq 2} (-1)^{|I|}
\nu( \mathop{\cap}_{\ell \in I} E_\ell)$. Note that
$$
|\phi(I')|
\leq 2^k \max_{1 \leq i < j \leq k} \nu\big( E_i \cap E_j \big).
$$
Applying the bound $\big| e^x - \sum_{n=0}^M \frac{x^n}{n!} \big| \leq
\frac{|x|^{M+1}}{(M+1)!} e^{|x|}$ to $e^{\phi(I')}$ in (\ref{corrbd}) then gives \eqref{eq:Poisvac.orderM}.
\qed
\bigskip

\medskip

\noindent {\bf Acknowledgements.} 
We would like to thank an anonymous reviewer for her/his 
careful reading and insightful comments which helped to clarify the presentation of 
some subtle points in the construction. M.B.\ would like to thank Anton Wakolbinger,
who originally posed the question addressed in this manuscript, see
\cite{B03}, for his encouragement and many interesting
discussions. M.B.\ would also like to thank Ted Cox and Jeremy Quastel
for stimulating discussions and ideas that after many years' gestation
took the form presented here.

R.S.\ is supported by AcRF Tier 1
grant R-146-000-185-112.  M.B.\ is in part supported by 
DFG priority programme SPP 1590 Probabilistic structures in evolution 
through grant BI 1058/3-1.

% %%%%%%%%%%%%%%%%%%%%%%%%%%% References %%%%%%%%%%%

\end{appendix}

\end{document}